\documentclass[10pt]{amsart}

\usepackage{mathtools}
\usepackage{amsmath}
\usepackage{amsfonts,amssymb}
\usepackage{enumerate}
\usepackage{mathrsfs}
\usepackage{amsthm}
\usepackage{color}
\usepackage{esint}

\newtheorem{theorem}{Theorem}

\newtheorem{definition}[theorem]{Definition}
\newtheorem{lemma}[theorem]{Lemma}
\newtheorem{proposition}[theorem]{Proposition}

\let\p=\partial

\let\O=\Omega

\newcommand{\be}{\begin{equation}}
\newcommand{\bm}{\begin{multline}}
\newcommand{\ee}{\end{equation}}
\newcommand{\dd}{\mathrm{d}}

\numberwithin{equation}{section}
\numberwithin{theorem}{section}

\def\p{\partial}

\def\O{\Omega}

\def\B{\begin{equation}}
\def\E{\end{equation}}
\def\BN{\begin{eqnarray*}}
\def\EN{\end{eqnarray*}}

\begin{document}
	
	
	\title{Initial value problem for the Free-boundary Magnetohydrodynamics with zero magnetic boundary condition}
	\author{Donghyun Lee}
	\maketitle
	
	\begin{abstract}
		We show local well-posedness of fluid-vacuum free-boundary magnetohydrodynamic(MHD) with both kinematic viscosity and magnetic diffusivity under the gravity force. We consider three-dimensional problem with finite depth and impose zero magnetic field condition on the free boundary and in vacuum. Sobolev-Slobodetskii space (Fractional Sobolev space) is used to perform energy estimates. Main difficulty is to control strong nonlinear couplings between velocity and magnetic fields. In \cite{DH}, we send both kinematic viscosity and magnetic diffusivity to zero with same speed to get ideal (inviscid) free-boundary magnetohydrodynamics using the result of this paper.
	\end{abstract}
	
	\section{Introduction}
	Magnetohydrodynamics (MHD) is a mathematical model for highly collisional charged particles. For example, plasma in tokamak, a device for nuclear fusion, is studied by MHD equations. In MHD model, the behavior of fluid particles is governed by Newton's law under the effect of electromagnetic force which is described by Maxwell' equations. Also, the magnetic field is also affected by electric field which is generated by movement of charged particles. In mathematical view point, MHD is described by two equations about velocity and magnetic fields. For velocity field, Navier-Stokes (or Euler) equation is considered under the effect of the electro magnetic Lorentz force. The evolution of magnetic field is governed by the Faraday's law with help of the Ohm's law and the Ampere's law. We note that displacement current in electromagnetism is ignored in MHD theory. In this paper, we consider MHD system with free boundary which has attracted a lot of attention for a long time. Free boundary problem describes the motion of fluid or plasma underneath vacuum or other fluid. Unlike to fixed boundary problem, we consider moving boundary problem. Therefore, the profile of the free boundary is also a function which solves some PDEs. Moreover, to determine the dynamics of the free boundary, we impose two more boundary conditions: kinematic boundary condition and continuity of stress tensor. Kinematic boundary condition is an evolution PDE of the boundary profile and means that a particle on the free surface stays on the surface as long as the solution is smooth. Continuity of stress tensor gives the balance of stress tensors on the interface of two different fluids. These two conditions will be explained again in this section.  \\
	
	\indent Before, we introduce problem setup, let us explain brief history of free-boundary Navier-Stokes and free-boundary MHD. Free-boundary Naiver-Stokes problem has a long story. In 1980, local regularity of free-boundary problem of Navier-Stokes equation was solved by Beale \cite{JB}, in the function space $K^r := L^2H^r \cap H^{ \frac{r}{2}}L^2$. Taking surface tension into account, A.Tani \cite{AT} proved local existence of the free-boundary problem. For global existence result, we introduce \cite{JBL} by Beale, \cite{ATNT} by A.Tani and N.Tanaka, and \cite{GuoTice} by Y.Guo and I.Tice. All these papers solved problems in unbounded  domains (at least horizontally unbounded domain). For finite mass case (isolated mass), V.A.Solonnikov solved local existence of the free-boundary problem in H\"older space in \cite{VS1}. There are many another articles by Solonnikov, Termamoto, and Allain. We refer introduction of \cite{AT} for more information, for both viscous cases and inviscid cases. Moreover, we also refer \cite{NMFR1} for additional references including irrotational inviscid free-boundary problem.   \\
	\indent For free-boundary magnetohydrodynamics (MHD), we should consider magnetic force effect for fluid equation. Moreover, we should treat one more time evolution PDE for magnetic field which looks like transport equation with strong couplings with velocity field. We will see that it looks like Navier-Stokes without pressure term. Physical properties of the conducting fluid determines whether the equation will be hyperbolic (perfect conducting fluid) or parabolic (not perfect-conducting fluid). This new equation will be derived from the Maxwell's equations of electromagnetic theory. Main difficulty of MHD equation is that we see high order nonlinear couplings between two vector fields. From this reason, considering diffusion effect both in velocity and magnetic field (non-perfect conducting viscous fluid) gives big advantage in solving the problem and many results were obtained in viscous cases unlike to free-boundary Euler. We refer \cite{PS2} for isolated mass setting when surface tension exists. They solved the problem using Hanzawa coordinate transform to the free boundary to use advantage of surface tension. We also refer apriori estimate result \cite{HAO} by Hao and Luo for ideal(inviscid) MHD case. In this paper, we solve free-boundary MHD using Lagrangian coordinate in finite depth domain without taking surface tension into account. On the free-boundary, we impose zero magnetic field condition. We refer \cite{DH} for an inviscid case result of free-boundary MHD which was obtained through vanishing viscosity technique.  \\
	
	\indent Let us formulate viscous free-boundary MHD. We use $u=(u_1,u_2,u_3)$ and $H=(H_1,H_2,H_3)$ to denote velocity and magnetic fields in Eulerian coordinate, respectively. Moving fluid domain at time $t$ is defined by fixed bottom and free boundary profile, 
	\[
		\O(t) := \{ (x,y,z)\in\mathbb{R}^{3} \ \vert \ -1 \leq z \leq h(t,x,y) \}\quad \text{and}\quad \O := \O(t=0),
	\] 
	where $h(t,x,y)$ is the profile of the free boundary. Note that the velocity $u$ is defined in $\O(t)$ and magnetic field $H$ is defined in both fluid domain $\O(t)$ and vacuum region. For vacuum domain, we use $\O_{V}(t)$ to denote
	\[
		\O_{V}(t) = \{ (x,y,z)\in\mathbb{R}^{3} \ h(t,x,y) < z \}\quad \text{and}\quad \O_{V} := \O_{V}(t=0) .
	\]
	We also define the free surface as 
	\[
	S_{F}(t) := \p\O(t) = \{ (x,y,z)\in\mathbb{R}^{3}  \ \vert \ z = h(t,x,y) \}\quad \text{and} \quad S_{F} := S_{F}(0) = \p\O,
	\] 
	and fixed bottom
	\[
		S_{B}(t) = S_{B} := \{ (x,y,z)\in\mathbb{R}^{3} \ \vert \ z=-1 \}.
	\]
	 We solve three-dimensional incompressible viscous magnetohydrodynamics(MHD) system with finite depth $z=-1$ and constant downward gravitational force,
	\begin{equation} \label{1.1}
	\begin{cases}
	\partial_t u + (u\cdot\nabla)u + \nabla P = \triangle u + (H\cdot\nabla)H - \frac{1}{2}\nabla|H|^2,\quad\text{in}\quad\Omega(t),  \\  
	\partial_t H + (u\cdot\nabla)H = \triangle H + (H\cdot\nabla)u  ,\quad\text{in}\quad\Omega(t), \\  
	\nabla\cdot u = 0, \quad \text{in}\quad \O(t),  \\
	\nabla\cdot H = 0,\quad\text{in}\quad \Omega(t) \cup S_{F}(t) \cup \O_{V}(t), \\
	P\mathbf{n} - 2\mathbf{D}(u)\mathbf{n}= gh\mathbf{n} + (H \otimes H - \frac{1}{2}I|H|^2)\mathbf{n},  \quad\text{on}\quad S_F(t), \\
	\partial_t h = u\cdot\mathbf{N},\quad\text{in}\quad S_{F}(t),  \\
	H = 0,\quad\text{on}\quad S_F(t) \cup \O_{V}(t), \\
	u = H = 0,\quad \text{on} \quad S_{B},  \\
	\end{cases}
	\end{equation}
	with initial data $u(0) = u_{0}, \ H(0) = H_{0}$, where $\mathbf{N}:= (-\partial_1 h, -\partial_2 h, 1)$ means outward normal vector on the free-boundary. We use $\mathbf{n}$ to denote the outward unit normal vector $\frac{\mathbf{N}}{|\mathbf{N}|}$ and $\mathbf{D}(u)$ to denote symmetric part of $\nabla u$,  
	\[
	\mathbf{D}(u) := \frac{\nabla u + (\nabla u)^T}{2}.  \\
	\]
	Let us explain physical meanings of the system (\ref{1.1}).  \\
	\textit{\textbf{1st} equation in (\ref{1.1})}: From classical electromagnetic theory, the Lorentz force is generated when a charged particle moves. The Lorentz force is described by divergence of magnetic stress tensor $T_{M}$,  
	\begin{equation}\label{mag tensor}
		T_{M} := (H\otimes H) - \frac{1}{2}|H|^{2}I,
	\end{equation}
	where $I$ is $3\times 3$ identity matrix. Note that $\nabla\cdot T_{M} = (H\cdot\nabla)H - \frac{1}{2}\nabla|H|^2$ with help of divergence-free property $\nabla\cdot H = 0$. Therefore, the first equation in (\ref{1.1}) describes Navier-Stokes equation under the effect of Lorentz force. We note that $P$ is the pressure combined with constant downward gravitational force $g$. This is why we see gravitation $g$ in 5th equation of (\ref{1.1}) in stead of the first equation.  \\
	
	\noindent \textit{\textbf{2nd} equation in (\ref{1.1})}: This is an evolution PDE for magnetic field $H$. From classical electromagnetic theory, this describes the Faraday's law,
	\begin{equation} \label{Faraday}
		\frac{\p H}{\p t} = -\nabla\times \mathbf{E},
	\end{equation}
	where $\mathbf{E}$ is electric field. Meanwhile, from the Ohm's law and the Ampere's law (without displacement current), we have
	\[
		\mathbf{J} = \sigma (\mathbf{E} + u\times H),\quad \nabla\times H = \mu\mathbf{J},
	\]
	where $\sigma$ is electric conductivity, $\mu$ is magnetic permeability, and $\mathbf{J}$ is current density. Therefore, from (\ref{Faraday}),
	\begin{equation} \label{Fara2}
	\begin{split}
		\frac{\p H}{\p t} &= -\nabla\times [ \frac{\mathbf{J}}{\sigma} - u\times H ] = \nabla\times [u\times H - \frac{1}{\mu\sigma} \nabla\times H]  \\
		&= \frac{1}{\mu\sigma} \Delta H + \nabla\times (u\times H)  \\
		&= \frac{1}{\mu\sigma} \Delta H + u(\nabla\cdot H) - H(\nabla\cdot u) + (H\cdot\nabla)u - (u\cdot\nabla)H \\
		&= - (u\cdot\nabla)H + (H\cdot\nabla)u + \frac{1}{\mu\sigma} \Delta H,
	\end{split}
	\end{equation}
	where we used vector identity $\nabla\times\nabla\times H = -\Delta H$ with divergence-free. We refer \cite{Dav} for more in details and derivations. Note that ideal MHD assumes perfect conducting fluid, $\sigma=\infty$, so we have no magnetic diffusivity term in (\ref{Fara2}) and this gives inviscid ideal MHD, e.g. \cite{HAO} and \cite{DH}.    \\
	
	\noindent \textit{\textbf{3rd} and \textbf{4th} equations in (\ref{1.1})}: Third equation means well-known incompressibility of the fluid. Fourth equation in (\ref{1.1}) comes from classical electromagnetic theory which implies that there is no magnetic monopole. \\
	
	\noindent \textit{\textbf{5th} equation in (\ref{1.1})}: Fifth equation in (\ref{1.1}) means continuity of stress tensor on the free boundary $S_{F}(t)$. Since fluid stress tensor is given by 
	\[
		T := 2\mathbf{D}u - PI, 	\\
	\]
	boundary condition $T\mathbf{n} + T_{M}\mathbf{n} = -gh\mathbf{n}$ yield 5th equation in (\ref{1.1}), where $T_{M}$ is defined in (\ref{mag tensor}).  \\
	
	\noindent \textit{\textbf{6th} equation in (\ref{1.1})}: Sixth equation of (\ref{1.1}) means kinematic boundary equation on the free surface. This is equivalent to 
	$$
	\frac{DF}{Dt} = 0\,\,\,\,\text{on}\,\,\,\,S_F(t),
	$$
	where $\frac{D}{Dt}$ is material derivative and $F(t,x,y,z) = 0$ for $z = h(t,x,y)$. Therefore, a particle on the free surface at initial time stays on the free surface $S_{F}(t)$ at later time $t$ also.  \\
	
	\noindent \textit{\textbf{7th} equation in (\ref{1.1})}: Let us use $H_{vac}$ and $H_{plasma}$ to denote magnetic fields in vacuum and plasma regions, respectively. Since displacement current is not assumed in MHD theory, magnetic field in vacuum solves
	\[
	\nabla\cdot H_{vac} = 0,\quad \nabla\times H_{vac} = 0,
	\]
	in general. Moreover, from parabolic property of the second equation in (\ref{1.1}), magnetic field $H$ propagates into vacuum region eventhough we assume zero magnetic field initially. Therefore, $H=0$ condition on the free surface and in the vacuum should be understood in the sense of imposed constraint, instead of natural propagation. In the real world, we may need some equipments to achieve this constraint. However, this is not harmful mathematically, as we impose zero boundary condition for heat equation on boundary for a fixed domain. Note that linearized equation of the second equation in (\ref{1.1}) in Lagrangian coordinate is simple heat equation. \\
	
	\indent Let us mention about continuity of the magnetic field on the free surface. We are considering nonzero magnetic diffusivity from $\Delta H$. This implies that charged fluid (or plasma) is not perfect conductor. Physically, this corresponds to zero surface current. This yields, 
	\[
		\nabla\times H = 0,\quad\text{on}\quad S_{F}(t),
	\]
	and applying curl-free property to a closed infinitesimal circuit near the boundary, we get tangential continuity of magnetic field $H$ on the free surface $S_{F}(t)$, 
	\begin{equation} \label{tan conti}
	\big(\mathbf{I-\mathbf{n}\otimes\mathbf{n}} \big)H_{vac} = \big(\mathbf{I-\mathbf{n}\otimes\mathbf{n}} \big) H_{{plasma}},\quad\text{on}\quad S_{F}(t).
	\end{equation}
	Similarly, using divergence-free property and considering a cylinder with infinitesimal height crossing the free boundary, we easily get normal continuity
	\begin{equation} \label{nor conti}
	H_{vac}\cdot\mathbf{n} = H_{{plasma}}\cdot\mathbf{n},\quad\text{on}\quad S_{F}(t).
	\end{equation}
	From (\ref{tan conti}) and (\ref{nor conti}), we get continuity of magnetic field $H$ acrossing the free surface $S_{F}(t)$. We note that we have only normal continuity (\ref{nor conti}) in the case of ideal MHD, which contains hyperbolic equation,
	\begin{equation} \label{hyp}
		\partial_t H + (u\cdot\nabla)H = (H\cdot\nabla)u .
	\end{equation}
	Tangential discontinuity comes from surface current which implies that $\O(t)$ is filled with perfect conducting fluid. For this case, we cannot give constraint $H=0$ in $\O_{V}(t)$, since PDE for $H$ is hyperbolic. In this case, geometric boundary condition $H\cdot\mathbf{n} = 0$ propagates from initial condition as far as we have a smooth solution for (\ref{hyp}). We refer \cite{HAO} for explanation and proof.  \\
	
	Therefore, we conclude that (\ref{1.1}) is proper approximate system for inviscid free-boundary MHD, when initial magnetic field is zero on the free boundary and vacuum. As explained in above paragraph, zero value of initial magnetic field on the free boundary propagates along the hyperbolic equation (\ref{hyp}) and there will be no magnetic flux into the vacuum region which means $H_{vac} = 0$ for all time. To obtain the solution of ideal (inviscid) MHD with zero initial magnetic field, uniform (in viscosity and diffusivity) regularity of (\ref{1.1}) and their limit was studied in \cite{DH}.  \\
	\indent We also give a remark on divergence-free condition on $H$. In the Navier-Stokes equation, pressure plays a role of Lagrange multiplier for constraint $\nabla\cdot u =0$. However, in Faraday's law, there is no such terms.  Instead, divergence-free condition of $H$ should be understood in the sense of propagation from initial compatibility condition. We should justify that our solution $H(t)$ satisfies $\nabla\cdot H = 0$ for $t < T$, where $[0,T)$ is time interval for local existence. Justification for divergence-free condition is given in the last section.  \\ 
	
	To simplify (\ref{1.1}), we define the total pressure $p$ as a sum of $P$ and magnetic pressure,  
	$$
	p := P + \frac{1}{2}|H|^2.
	$$
	Then (\ref{1.1}) is written by 
	\begin{equation} \label{system}
	\begin{cases}
	\partial_t u + (u\cdot\nabla)u - (H\cdot\nabla)H + \nabla p = \triangle u,\quad \text{in}\quad\Omega(t), \\
	\partial_t H + (u\cdot\nabla)H - (H\cdot\nabla)u  = \triangle H,\quad\text{in}\quad\Omega(t), \\
	\nabla\cdot u = 0, \quad \text{in}\quad \O(t),  \\
	\nabla\cdot H = 0,\quad\text{in}\quad \Omega(t) \cup S_{F}(t) \cup \O_{V}(t), \\
	p\mathbf{n} - 2\mathbf{D}(u)\mathbf{n}= gh\mathbf{n},\quad\text{on}\quad S_F(t), \\
	\p_{t} h = u^b\cdot\mathbf{N},\quad\text{on}\quad S_{F}(t), \\
	H = 0,\quad\text{on}\quad S_F(t) \cup \O_{V}(t), \\
	u = H = 0,\quad \text{on} \quad S_{B}.
	\end{cases}
	\end{equation}
	Initial data $u(0) = u_{0}$ and $H(0) = H_{0}$ must be given so that satisfy the following initial compatibility conditions
	\begin{equation} \label{compatibility}
	\begin{cases}
	\nabla\cdot u_0 = 0, \quad\text{in} \quad \Omega, \\
	\nabla\cdot H_0 = 0, \quad\text{in} \quad \Omega, \\
	\mathbf{\Pi}\mathbf{D}(u_0)\mathbf{n} = 0,\quad\text{on} \quad S_{F}, \\
	H_0 = 0,\quad\text{on} \ S_F \cup \O_{V}, \\
	u_{0} = H_{0} = 0,\quad \text{on} \quad S_{B},
	\end{cases}
	\end{equation}
	where $\mathbf{\Pi}$ means tangential projection operator $\mathbf{\Pi} := I - \mathbf{n}\otimes\mathbf{n}$. \\ 
	
	\subsection{Lagrangian coordinate and functional framework}
	To solve the problem in the fixed domain, we transform the problem (\ref{system}) into the initial domain using Lagrangian map $X(t,\cdot):\xi \rightarrow x$ which solves $\p_{t}X = u$. We introduce new variables $v$, $B$, and $q$ for velocity, magnetic field, and total pressure in Lagrangian coordinate,  
	\begin{equation} \label{lagrangian}
	\begin{split}
	v(t,\xi) &:= u(t,X(t,\xi)) = u(t,x),  \\
	B(t,\xi) &:= H(t,X(t,\xi)) = H(t,x),  \\
	q(t,\xi) &:= p(t,X(t,\xi)) = p(t,x).
	\end{split}
	\end{equation}
	
	\noindent Using notation $\p_{\xi_{i}} := \p_{i}$, we define
	\begin{equation} \label{loc deriva}
	\begin{split}
		\begin{bmatrix}
			\p_{x_{1}} \\ \p_{x_{2}} \\ \p_{x_{3}}  
		\end{bmatrix} 
		=
		\begin{bmatrix}
			\  &  \  &  \  \\  
			\  &  \frac{\p X}{\p \xi}  &  \  \\  
			\  &  \  &  \  \\  
		\end{bmatrix}^{-T} 
		\begin{bmatrix}
			\p_{1} \\ \p_{2} \\ \p_{3}  
		\end{bmatrix}
		:=
		\begin{bmatrix}
		\  &  \  &  \  \\  
		\  &  \mathcal{G}  &  \  \\  
		\  &  \  &  \  \\  
		\end{bmatrix}\
		\begin{bmatrix}
		\p_{1} \\ \p_{2} \\ \p_{3}  
		\end{bmatrix},
		\quad
		\p_{x_{i}} := \p_{v,i} := \sum_{j=1}^{3} {G}_{ij}\p_{j},
	\end{split}
	\end{equation}
	where $G_{ij}$ is $(i,j)$ entry of the matrix $\mathcal{G}$.  \\

	\noindent We should rewrite the system (\ref{system}) in terms of $(v,B,q)$ and $(t,\xi)$. Let $\Pi^x_{\xi}$ be corresponding transform from $(\Omega(t), \ S_F(t))$ to $(\Omega, \ S_F)$, then (\ref{system}) can be rewritten in terms of $(t,\xi)$ in $\Omega$. 
	\begin{equation} \label{Larg system}
	\begin{cases}
	v_t - \triangle_v v + \nabla_v q = B\cdot\nabla_v B\,\quad\text{in}\quad \O, \\
	B_t - \triangle_v B = B\cdot\nabla_v v,\quad\text{in}\quad \O, \\
	\nabla_v\cdot v = 0,\quad\text{in}\quad \O, \\
	\nabla_v\cdot B = 0,\quad\text{in}\quad \O\cup S_{F} \cup \O_{V}, \\
	q - 2\mathbf{D}_v(v)\mathbf{n}^{(v)}\cdot\mathbf{n}^{(v)} = gh\quad \text{on}\quad S_{F}, \\
	B = 0,\quad \text{on}\quad S_{F}\cup \O_{V}	\\
	v= H = 0,\quad\text{on}\quad S_{B},
	\end{cases}
	\end{equation}
	with initial compatibility conditions
	\begin{equation} \label{Larg compatibility}
	\begin{cases}
	\nabla_v \cdot v_0 = 0, \quad\text{in} \quad \O, \\
	\nabla_v \cdot B_0 = 0, \quad\text{in} \quad \O, \\
	\mathbf{\Pi}^{(v)}\mathbf{D}_v(v_0)\mathbf{n}^{(v)} = 0,\quad\text{on} \quad S_F, \\
	B_0 = 0,\quad\text{on} \quad S_{F}\cup \O_{V}, \\
	v_{0} = B_{0} = 0,\quad\text{on}\quad S_{B},
	\end{cases}
	\end{equation}
	where we used the notations
	\begin{equation} \label{symbols def}
	\begin{split}
	f^{(v)} &= f^{(v)}(t,\xi) = \Pi_{\xi}^x f(t,x),  \\
	\nabla_v &:= (\p_{v,1}, \p_{v,2}, \p_{v,3}) ,   \\
	\mathbf{D}_v(f) &:= \frac{1}{2} \Big( (\nabla_v f) + (\nabla_v f)^{T} \Big).  \\
	\end{split}
	\end{equation}
	Note that $\p_{v,i}$ is defined in (\ref{loc deriva}). We do not specify coordinate transform in vacuum region $\O_{V}(t)$, since $B(t)$ is constrained to be zero in $\O_{V}(t)$. We suffice to consider vector field in fluid domain $\O$.  	\\

	Before we explain main scheme of the proof, we introduce Sobolev-Slobodetskii space which justifies fractional derivatives eventhough domain is not whole $\mathbb{R}^{3}$. Note that this space is equivalent to standard $H^{s}(\mathbb{R}^{3})$ which is defined by Fourier transform, see \cite{HJR}. Throughout this paper we will use space-time domain 
	\begin{equation} \label{st domain}
		Q_{T} := \O \times [0,T) \quad\text{and}\quad S_{F,T} := S_F\times [0,T).
	\end{equation}
	
	\begin{definition} \label{basic W}
		By $W_2^l(\Omega)$, we define,	
		\begin{equation} \label{slobodetskii}
		\|u\|^{2}_{W_{2}^{l}(\Omega)} := \sum_{0 \leq |\alpha| < {l}}\|D^{\alpha} u\|^{2}_{L^{2}(\Omega)} + \|u\|^{2}_{\dot{W}_{2}^{l}(\Omega)},
		\end{equation}
		where
		$$
		\|u\|^{2}_{\dot{W}_{2}^{l}(\Omega)} := 
		\begin{dcases*}
		\sum_{|\alpha|=l}\|D^{\alpha} u\|^{2}_{L^2(\Omega)}\quad  \text{if} \quad l\in \mathbb{Z}, \\
		\sum_{|\alpha|=[l]} \int_{\Omega}\int_{\Omega} \frac{|D^{\alpha} u(x) - D^{\alpha} u(y)|^{2}}{ |x-y|^{ n+2\{l\} } }\quad  \text{if} \quad l= [l] + \{l\} \notin \mathbb{Z},\quad 0< \{l\} <1 .
		\end{dcases*}
		$$
	\end{definition}
	
	\begin{definition} \label{aniso W}
		We also define the anisotropic space $W_{2}^{l, \frac{l}{2} }(Q_T)$ as 
		\begin{equation} \label{aniso}
		\begin{split}
		\|u\|^{2}_{W_{2}^{l,\frac{l}{2}}(Q_T)} & :=  \|u\|^{2}_{W_{2}^{l,0}(Q_T)} + \|u\|^{2}_{W_{2}^{0,\frac{l}{2}}(Q_T)} \\
		& := \int_{0}^{T} \|u(t,\cdot)\|^{2}_{W_{2}^{l}(\Omega)} dt + \int_{\Omega} \|u(\cdot,x)\|^{2}_{W_{2}^{\frac{l}{2}}(0,T)} dx.
		\end{split}
		\end{equation}
	\end{definition}
	
	\begin{definition} \label{decay H}
		By $H^{l,\frac{l}{2}}_\gamma(Q_{T}),\gamma \geq 0$, we mean
		\begin{equation} \label{H space}
		\|u\|^{2}_{H^{l,\frac{l}{2}}_\gamma(Q_{T})} := \|u\|^{2}_{H^{l,0}_\gamma(Q_{T})} + \|u\|^{2}_{H^{0,\frac{l}{2}}_\gamma(Q_{T})},
		\end{equation}
		where
		\begin{equation}
		\begin{split}
		\|u\|^2_{H^{l,0}_\gamma(Q_T)} &:= \int_0^T e^{-2\gamma t} \|u(t,\cdot)\|^2_{\dot{W}_2^l(\Omega)}dt,  \\
		\|u\|^2_{H^{0,\frac{l}{2}}_\gamma(Q_T)} &:= \gamma^l\int_0^T e^{-2\gamma t} \|u(t,\cdot)\|^2_{L^2(\Omega)}dt  \\
		&+ \int_0^T e^{-2\gamma t} dt \int_0^\infty \left\| \left( \frac{\partial}{\partial t}\right)^k u_0(t,\cdot) - \left( \frac{\partial}{\partial t}\right)^k u_0(t-\tau,\cdot) \right\|^2_{L^2(\Omega)}\frac{d\tau}{\tau^{1+l-2k}},
		\end{split}
		\end{equation}
		if $\frac{l}{2}$ is not an integer, $k=[\frac{l}{2}], \ u_0(t,x) = u(t,x)(t>0)$, and $\ u_0(t,x) = 0(t<0)$. If $\frac{l}{2}$ is an integer, then the double integral in the norm should be replaced by 
		$$
		\int_{-\infty}^T e^{-2\gamma t} \left\| \left( \frac{\partial}{\partial t} \right)^{\frac{l}{2}}u(t,\cdot) \right\|^2_{L^2(\Omega)} dt,
		$$
		and $\left( \frac{\partial}{\partial t}\right)^j u |_{t=0} = 0 \ (j=0,1,2,\cdots,\frac{l}{2}-1)$ should be satisfied.
	\end{definition}
	
	\begin{definition}
		We also introduce the space $H_\gamma^{\ell+\frac{1}{2},\frac{1}{2},\frac{\ell}{2}}(S_{F,T})$ to treat trace on the boundary $S_{F,T}$.
		\begin{equation} \label{trace space}
		\begin{split}
		\|u\|^2_{H^{\ell+\frac{1}{2},\frac{1}{2},\frac{\ell}{2}}_\gamma(S_{F,T})} &:= \int_0^T e^{-2\gamma t} \left( \|u\|^2_{\dot{W}_2^{\ell+\frac{1}{2}}(S_{F})} + \gamma^{\ell} \|u\|^2_{\dot{W}_2^{ \frac{1}{2}}(S_{F})} \right) dt   \\
		& + \int_0^T e^{-2\gamma t} dt \int_0^\infty \left\| \left( \frac{\partial}{\partial t}\right)^k u_0(t,\cdot) - \left( \frac{\partial}{\partial t}\right)^k u_0(t-\tau,\cdot) \right\|^2_{W_2^{\frac{1}{2}}(S_{F})}\frac{d\tau}{\tau^{1+\ell-2k}},
		\end{split}
		\end{equation}
		if $\frac{\ell}{2}$ is not an integer, $k=[\frac{\ell}{2}], \ u_0(t,x) = u(t,x)(t>0)$, and  $u_0(t,x) = 0(t<0)$. If $\frac{\ell}{2}$ is an integer, then the double integral in the norm should be replaced by 
		$$
		\int_{-\infty}^T e^{-2\gamma t} \left\| \left( \frac{\partial}{\partial t} \right)^{\frac{\ell}{2}}u(t,\cdot) \right\|^2_{W_2^{\frac{1}{2}}(S_{F})} dt,
		$$
		and $\left( \frac{\partial}{\partial t}\right)^j u |_{t=0} = 0 \ (j=0,1,2,\cdots,\frac{\ell}{2}-1)$ should be satisfied. 
	\end{definition}
	
	\begin{definition} \label{C space}
		For $l \in(\frac{1}{2},1)$, we define
		\begin{equation} \label{mixed space}
		(\|u\|_{Q_T}^{(l+2)})^{2} := (\|\p_{t}u\|_{Q_T}^{(l)})^{2} + \sum_{|s|=2} (\|D_x^s u\|_{Q_T}^{(l)})^{2} + \sum_{|s|=0}^1\|D_x^s u\|_{L^2(Q_T)}^{2} ,
		\end{equation}
		where 
		\begin{equation} \label{QT space} 
		(\|u\|_{Q_T}^{(l)})^{2} := \|u\|^2_{W_2^{l,\frac{l}{2}}(Q_T)} + T^{-l}\|u\|^2_{L^2(Q_T)}.
		\end{equation}
		And, since we treat nonlinear terms on right hand side of (\ref{Larg system}), we need $L^{\infty}_{T}$ type norm. We define
		\begin{equation} \label{infty space}
		\|u\|^2_{H^{l+2,\frac{l}{2}+1}(Q_T)} := \big( \|u\|^{(l+2)}_{Q_T}\big)^{2} + \sup_{t<T}\|u(t)\|^2_{W_2^{l+1}(\Omega)}.
		\end{equation}
	\end{definition}
	

	\subsection{Scheme of proof and main result} 
	In section 2, we state some preliminary lemmas for basic estimates of Lagrangian map and algebraic properties of Sobolev-Slobodetskii space. From the first equation in (\ref{Larg system}), we should solve linear Stokes problem which is obtained by replacing all differential operators $\nabla_v$ and $\Delta_v$ into $\nabla$ and $\Delta$. At the same time, we just get simple heat equation with zero boundary condition from the second equation in (\ref{Larg system}). Linear Stokes equation is well-developed already as we can see in \cite{JB} and \cite{AT}. In section 3, we modify these theories to obtain linear solvability in the space $H^{l+2,\frac{l}{2}+1}$. Unlike to \cite{JB} and \cite{AT}, $H^{l+2,\frac{l}{2}+1}$ contains $L^{\infty}_{T}$ type energy term and this plays a critical role in nonlinear estimates which look like $B\cdot\nabla B$ and $B\cdot\nabla v$ on the RHS of (\ref{Larg system}).  \\
	\indent In section 4, we obtain main nonlinear solvability of Stokes equation and Heat equation. If we have only $L^{2}_{T}$ type regularity without $L^{\infty}_{T}$ terms, then we cannot control nonlinear terms because 
	\[
		\int_{0}^{T} \| B\cdot\nabla B \|_{H^{\ell}} \sim  \|B\|^{2}_{L^{2}_{T} H^{\ell+1}},
	\]
	which does not include small factor $T^{s}$ for some $s > 0$. Instead, with help of $L^{\infty}_{T}$ type regularity, the function space with fractional order and algebraic properties in Lemma \ref{interpolation lemma} play critical roles in deriving small factor $T^{s > 0}$ in nonlinear estimate, 
	\[
		\int_{0}^{T} \| B\cdot\nabla B \|_{H^{\ell}} \sim  T^{s} \|B\|^{2}_{L^{2}_{T} H^{\ell+2}\cap L^{\infty}_{T}H^{\ell}},\quad\text{for some} \quad s > 0.
	\]
	Therefore we can perform standard contraction mapping argument to get nonlinear solution. Sharp nonlinear estimate with Definitions \ref{basic W} -- \ref{C space} is proved in Lemma \ref{nonlinear est}.  \\
	\indent In section 5.1, we solve fully nonlinear problem (\ref{fully nonlin system}) from solution of (\ref{cst coeff nonlin}). We change all the differential operators $\nabla$ and $\Delta$, into $\nabla_v$ and $\Delta_v$. As introduced in Lemma \ref{nonlinear est}, we also need similar sharp estimate Lemma \ref{fully nonlin est} to apply contraction mapping principle. \\
	\indent Meanwhile, divergence-free condition of magnetic field is not considered when we solve heat equation, because heat equation does not contain pressure term, which is Lagrange multiplier for divergence-free condition. Therefore, we should claim that divergence-free condition propagates from initial time, if $\nabla\cdot H_0 = 0$ at time $t=0$. We can find a very good nonlinear cancellation between $(u\cdot\nabla)H$ and $(H\cdot\nabla)u$, and derive a convection-diffusion equation for divergence of $H$. Maximum principle of convection-diffusion equation yields divergence-free property of $H$ directly. This is explained in the section 5.2.  \\
	
	Now we state the main result of this paper.
	\begin{theorem} \label{main theorem}
		Let $l \in (\frac{1}{2},1)$. For initial data $h_0 \in W_2^{l+5/2}(S_F), \ v_0 \in W_2^{l+1}(\Omega)$, and $\ B_0 \in W_2^{l+1}(\Omega)$ with initial compatibility conditions
		\begin{equation*}
		\begin{cases}
		\nabla\cdot v_0 = \nabla\cdot B_0 = 0,\quad\text{in}\quad\Omega, \\
		\mathbf{D}(v_0)\mathbf{n} - (\mathbf{D}(v_0)\mathbf{n}\cdot \mathbf{n})\mathbf{n} = 0,\,\,\,\,\text{on}\,\,\,\,S_{F}, \\
		B_0 = 0,\quad\text{on}\quad S_{F},  \\
		v_{0} = B_{0} = 0,\quad\text{on}\quad S_{B},
		\end{cases}
		\end{equation*}
		there exist a unique solution 
		\[
			v \in H^{l+2, \frac{l}{2}+1}(Q_{T^*}),\quad B \in H^{l+2, \frac{l}{2}+1}(Q_{T^*}),\quad \nabla q \in W_2^{l, \frac{l}{2}}(Q_{T^*}),\quad q \in W_2^{l+ \frac{1}{2}, \frac{l}{2}+ \frac{1}{4}}(S_{F,{T^*}}),
		\]
		to the system (\ref{Larg system}) for some $T^{*} > 0$. Moreover, we have the following estimate,
		\begin{equation}
		\begin{split}
		& \|v\|_{H^{l+2, \frac{l}{2}+1}(Q_{T^*})} + \|B\|_{H^{l+2, \frac{l}{2}+1}(Q_{T^*})} + \|\nabla q\|_{W_2^{l, \frac{l}{2}}(Q_{T^*})} + \|q\|_{W_2^{l+ \frac{1}{2}, \frac{l}{2}+ \frac{1}{4}}(S_{F,{T^*}})}  \\
		&\quad \leq C_0 \Big( \|v_0\|_{W_2^{l+1}(\Omega)} + \|B_0\|_{W_2^{l+1}(\Omega)} + \|h_0\|_{W_2^{l+\frac{3}{2}}(\mathbb{R}^2)} \Big).
		\end{split}
		\end{equation}
	\end{theorem}
	
	\noindent\textbf{Notation}
	\ Throughout this paper, $C$ is a positive constant and $C(\cdot,\cdot)$ is a positive constant depending \textit{increasingly} on its arguments. They may vary line to line.
	
	\section{Preliminaries}
	In this section, we state some preliminary lemmas for Lagrangian map, unit normal vector, and algebraic properties of (\ref{slobodetskii}) . Let us assume that $T$ satisfies
	\begin{equation} \label{T assumption}
	T^{ \frac{1}{2}} \|u\|_{Q_{T}}^{(l+2)} \leq \delta.
	\end{equation}
	
	The first lemma is standard for Lagrangian coordinate. $\mathcal{G} = \Big( \frac{\p X}{\p \xi} \Big)^{-T}$ is defined in (\ref{symbols def}) and we use $\mathcal{G}^{(u)}$ and $\mathcal{G}^{(u^{\prime})}$ to denote that Lagrangians with respect to velocity field $u$ and $u^{\prime}$ respectively.
	\begin{lemma} \label{larg lemma}
		When $u,u' \in W_2^{l+2, \frac{l}{2}+1}(Q_{T})$ satisfy (\ref{T assumption}), we have the following estimates for $t\leq T$: 
		\begin{equation*}
		\begin{split}
		\left\| \mathcal{G}^{(u)} - \mathcal{G}^{(u^{\prime})} \right\|_{W_{2}^{ l+1,\frac{l+1}{2} }(Q_t)} &\leq C(t,\delta)t^{1-\frac{l}{2}} \|u-u'\|_{W_2^{ l+2, \frac{l+1}{2} }(Q_t)},  \\
		\left\| \frac{\partial}{\partial t}\left(\mathcal{G}^{(u)} - \mathcal{G}^{(u')} \right) \right\|_{W_2^{l+1,\frac{l+1}{2}}(Q_t)} &\leq C(t,\delta) \|u-u'\|_{W_2^{l+2,\frac{l}{2}+1}(Q_t)},  \\
		\left\| \frac{\partial}{\partial t}\left(\mathcal{G}^{(u)} - \mathcal{G}^{(u')}\right) \right\|_{L^2(Q_t)} &\leq C(t,\delta)t^{\frac{1}{2}} \|u-u'\|_{W_2^{l+2,\frac{l}{2}+1}(Q_t)},
		\end{split}
		\end{equation*}
		where $C(t,\delta)$ is a positive constant depending increasingly on both arguments.
	\end{lemma}
	\begin{proof}
		Lemma 4.1 in \cite{AT}.
	\end{proof}
	
	We also need estimate about outward unit normal vector in Lagrangian framework.
	\begin{lemma} \label{normal lemma}
		Let $\delta_0$ be a positive root of the equation $6\delta^2 + 6\delta - 1 = 0$. If $u,u' \in W_2^{l+2,\frac{l}{2}+1}(Q_{T_0})$ satisfying (\ref{T assumption}) with $\delta < \delta_0$, and if $0<t\leq T_0$, then 
		\begin{equation}
		\begin{split}
		\|\mathbf{n}^{(u)} - \mathbf{n}^{(u^{\prime})}\|_{W_2^{l+\frac{1}{2},\frac{l}{2}+\frac{1}{4}}(S_{F,t})} &\leq C(t,\delta)t^{1-\frac{l}{2}} \|u - u^{\prime}\|_{W_2^{l+2,\frac{l}{2}+1}(Q_t)}, \\
		\left\|\frac{\partial}{\partial t}(\mathbf{n}^{(u)} - \mathbf{n}^{(u^{\prime})})\right\|_{W_2^{l+\frac{1}{2},\frac{l}{2}+\frac{1}{4}}(S_{F,t})} &\leq C(t,\delta) \|u - u^{\prime}\|_{W_2^{l+2,\frac{l}{2}+1}(Q_t)},
		\end{split}
		\end{equation}
		where $C(t,\delta)$ is a positive constant depending increasingly on both arguments.
	\end{lemma}
	\begin{proof}
		Lemma 4.3 in \cite{AT}.
	\end{proof}
	
	For Sobolev-Slobodetskii space, we have the following algebraic properties.
	\begin{lemma} \label{interpolation lemma}
		For smooth functions $u(x)$ and $v(x)$, defined in a domain $\Omega \subset \mathbb{R}^n$, they satisfy the following inequalities,
		\begin{equation}
		\begin{split}
		\|uv\|_{W_2^l(\Omega)} &\leq c \|u\|_{W_2^l(\Omega)} \|v\|_{W_2^s(\Omega)},\quad s > \frac{n}{2}, \quad l < \frac{n}{2},  \\
		\|uv\|_{L^2(\Omega)} &\leq c \|u\|_{W_2^l(\Omega)} \|v\|_{W_2^{n/2-l}(\Omega)},\quad l < \frac{n}{2}, \\
		\|uv\|_{W_2^l(\Omega)} &\leq c \left( \|u\|_{W_2^l(\Omega)} \|v\|_{W_2^s(\Omega)} + \|v\|_{W_2^l(\Omega)} \|u\|_{W_2^s(\Omega)}\right),\quad l,s > \frac{n}{2}.
		\end{split}
		\end{equation}
	\end{lemma}
	\begin{proof}
		We refer \cite{PS}.
	\end{proof}
	
	The following is Lemma 6.3 and its corollary in \cite{solo}.
	\begin{lemma} \label{solo lemma 63}
		\noindent (1) Let $r<1$ and $u\in W^{0, \frac{r}{2}}_{2}(Q_{T})$. Then the inequality
		\begin{equation}
		\begin{split}
		\int_{0}^{T} \|u(\cdot,t)\|^{2}_{L^{2}(\O)} \frac{1}{t^{r}} \dd t &\leq C_{1} T^{-r} \int_{0}^{T} \|u(\cdot,t)\|^{2}_{L^{2}(\O)}  \dd t  \\
		&+ C_{2} \int_{0}^{T} \int_{0}^{T} \|u(\cdot,t-\tau) - u(\cdot,t)\|^{2}_{L^{2}(\O)} \frac{1}{\tau^{1+r}} \dd t \dd\tau
		\end{split}
		\end{equation}
		holds. Above inequality holds also for $1<r<2$, except for the first term in the RHS is $u(0,x) = 0$. Two constants $C_{1}$ and $C_{2}$ are independent of time $T$ and velocity $u$.  \\
		\noindent (2) Let $r>1$ not be an odd integer. For any function $u\in W^{r,\frac{r}{2}}_{2}(Q_{T})$ satisfying the relation $\frac{\dd^{j}}{\dd t^{j}} u \vert_{t=0} \ (j=0,1,\cdots,[(r-1)/2])$, the following inequality holds:
		\begin{equation}
		\begin{split}
		\|u\|^{2}_{H^{r, \frac{r}{2}}_{0}(Q_{T})} &\leq C_{3} \Big\{ \int_{0}^{T} \|u(\cdot,t)\|^{2}_{\dot{W}^{r}_{2}(\O)} \dd t      \\
		&+ \int_{\O} \|u(\cdot,t)\|^{2}_{\dot{W}^{\frac{r}{2}}_{2}(0,T)} \dd x + T^{-r+2k} \Big\| \frac{\dd^{k}}{\dd t^{k}} u\Big\|^{2}_{L^{2}(Q_{T})}  \Big\},
		\end{split}
		\end{equation}
		where $k=[ \frac{r}{2}]$ and $C_{3}$ is time-independent.
	\end{lemma}

	\section{linear Stokes problem}
	From (\ref{Larg system}), we start from linear Stokes problem, 
	\begin{equation} \label{stokes system}
	\begin{cases}
	v_t - \triangle v + \nabla q = f, \quad \text{in}\quad Q_{T}, \\
	\nabla\cdot v = \rho,\quad \text{in}\quad Q_{T}, \\
	v(0) = v_0,\quad \text{in}\quad \Omega\times\{t=0\}, \\
	2 \mathbf{\Pi}\mathbf{D}(v) = d,\quad \text{on}\quad S_{F,T}, \\
	-q + 2\mathbf{D}(v)\mathbf{n}\cdot \mathbf{n} = b,\quad \text{on}\quad S_{F,T}, \\
	u = 0,\quad \text{on}\quad S_{B}.
	\end{cases}
	\end{equation}
	
	Since we are assuming flat bottom, we can apply Theorem 1 in \cite{AT} and state
	\begin{proposition} \label{stokes with zero initial}
		(Linear Stokes problem with zero initial data)
		Let $l>\frac{1}{2}$, $0<T\leq\infty$, and $\gamma$ be sufficiently large so that $\gamma\geq\gamma_0\geq 1$. And we assume $S_F \in W_2^{l+\frac{3}{2}}$. Then, for
		$$
		(f,\rho,d,b)\in H_{\gamma}^{l,\frac{l}{2}}(Q_T) \times H_{\gamma}^{l+1,\frac{l+1}{2}}(Q_T) \times H_{\gamma}^{l+\frac{1}{2},\frac{l}{2}+\frac{1}{4}}(S_{F,T}) \times H_{\gamma}^{l+\frac{1}{2}, \frac{1}{2}, \frac{l}{2}}(S_{F,T}),
		$$
		satisfying compatibility conditions $d\cdot \mathbf{n} = 0$, $(\rho,d)|_{t=0} = 0$, $\rho=\nabla\cdot R$, and $R\in H_{\gamma}^{l+1,1,\frac{l}{2}}(Q_T)$, there exist a unique solution $(v,q)\in H_{\gamma}^{l+2,\frac{l}{2}+1}(Q_T) \times H_{\gamma}^{l+1,1,\frac{l}{2}}(Q_T)$ to the problem (\ref{stokes system}) with zero initial condition $v_0 = 0$. Moreover we have the following estimate.
		\begin{equation*}
		\begin{split}
		\|v\|_{H_{\gamma}^{l+2,\frac{l}{2}+1}(Q_T)} + \|q\|_{H_{\gamma}^{l+1,1,\frac{l}{2}}(Q_T)}   
		& \leq  C \big\{ \|f\|_{H_{\gamma}^{l,\frac{l}{2}}(Q_T)} + \|\rho\|_{H_{\gamma}^{l+1,\frac{l+1}{2}}(Q_T)} + \|R\|_{H_{\gamma}^{0,\frac{l}{2}+1}(Q_T)} \\
		&\quad + \|d\|_{H_{\gamma}^{l+\frac{1}{2},\frac{l}{2}+\frac{1}{4}}(S_{F,T})} + \|b\|_{H_{\gamma}^{l+\frac{1}{2},\frac{1}{2},\frac{l}{2}}(S_{F,T})} \big\},
		\end{split}
		\end{equation*}
		where $\|q\|^2_{H_{\gamma}^{l+1,1,\frac{l}{2}}(Q_T)} := \|q\|^2_{H_{\gamma}^{l,\frac{l}{2}}(Q_T)} + \|\nabla q\|^2_{H_{\gamma}^{l,\frac{l}{2}}(Q_T)}$.
	\end{proposition}
	
	To solve the problem with general data $v_0 \neq 0$, we produce a function $U_0(t,x)$ which satisfies initial data $U_0(0) = v_0$ and their some space-time sobolev norms are bounded by $v_0$.
	
	\begin{lemma} \label{extention lemma}
		(1) For $v\in W_2^{l+2,\frac{l}{2}+1}(Q_T), \ l\in(\frac{l}{2},1)$, there is an extension $U \in W_2^{l+2,\frac{l}{2}+1}(Q_{\infty})$ such that
		$$
		\|U\|^{(l+2)}_{Q_\infty} \leq C \|v_0\|_{W_2^{l+1}(\Omega)},
		$$
		where $\|U\|_{Q_{\infty}}^{(l+2)}$ is defined in Definition (\ref{mixed space}).  \\
		(2) Let $w(0)=0$ and $w \in W_2^{l+2,\frac{l}{2}+1}(Q_T)$. Then there exist an extension $w_{\text{ext}} \in W_2^{l+2,\frac{l}{2}+1}(Q_\infty)$ such that
		$$
		\|w_{\text{ext}}\|_{W_2^{l+2,\frac{l}{2}+1}(Q_\infty)} \leq C \|w\|_{W_2^{l+2,\frac{l}{2}+1}(Q_T)}.
		$$
	\end{lemma}
	\begin{proof}
		By Lemma \ref{solo lemma 63} and trace theorem, there exist a vector field $U \in W_2^{l+2,\frac{l}{2}+1}(Q_{\infty})$ which satisfies $U(0) = v_0$ and the following inequality,
		\begin{equation} \label{U ineq}
		\int_{0}^\infty \Big( \|\partial_t U(\cdot,t)\|^2_{L^2(\Omega)} + \sum_{|s|=2} \|D^s U(\cdot,t)\|^2_{L^2(\Omega)}\Big)\frac{1}{|t|^l} dt + \|U\|^2_{W_2^{l+2,\frac{l}{2}+1}(Q_\infty)} \leq C \|v_0\|^2_{W_2^{l+1}(\Omega)}.
		\end{equation} 
		Since $U$ is not a function on a finite time interval, we can use interpolation with time-independent constant to get
		$$
		\|\partial_t U\|^2_{W_2^{l,0}(Q_{\infty})},\ \|D^2 U\|^2_{W_2^{0, \frac{l}{2}}(Q_{\infty})} \leq C \|U\|^2_{W_2^{l+2,\frac{l}{2}+1}(Q_{\infty})}.
		$$
		Hence, obviously,
		\begin{equation*}
		\|U\|^{(l+2)}_{Q_T} \leq C \|v_0\|_{W_2^{l+1}},\,\,\,\,\forall T>0.
		\end{equation*}
		For second one, we consider $w\in W_{2}^{s}(0,T,X)$ for some Hilbert space $X$ with integer $s$. We extend $w$ into $W_{2}^{l}(-\infty, T; X)$ by
		\begin{equation}
		w_{ext}(t) = 
		\begin{cases}
			0,\quad t<0,\\
			w(t),\quad 0\leq t < T,  \\
			3w(2T-t) - 2w(3T-2t),\quad t > T.
		\end{cases}
		\end{equation} 
		For non-integer $s$, we use interpolation to finish the proof. Applying this extension with $(s=0, \ X=W_{2}^{l+2})$ and $(s=\frac{l}{2}+1, \ X=L^{2})$, we get the result.
	\end{proof}

	Since $U$ is global in time function, so constant $C$ is time-independent. Using the function $U$, we solve the same linear problem with general initial data. We also need the estimate of $(v,B)$ in $C([0,T];W_2^{s})$ type norm, since nonlinear terms on the right hand side, $B\cdot\nabla B$ and $B\cdot\nabla v$ should be estimated in the form of $CT^{\gamma}\Lambda(\|v\|,\|B\|)$ (for some positive $\gamma>0$), where $\Lambda$ is some nice function. 

	\begin{proposition} \label{prop general stokes}
		(Linear stokes problem with general initial data)
		Let $l \in (\frac{1}{2},1)$, $0<T<\infty$, and $S_F \in W_2^{l+\frac{3}{2}}$. And $(f,\rho,u_0,(b,d))$ in (\ref{stokes system}) satisfies
		\begin{equation*}
		\begin{split}
		(f,\rho,u_0,(b,d)) &\in W_2^{l, \frac{l}{2}}(Q_T) \times W_2^{l+1, \frac{l+1}{2}}(Q_T) \times W_2^{l+1}(\Omega) \times W_2^{l+ \frac{1}{2}, \frac{l}{2}+ \frac{1}{4}}(S_{F,T}),  \\
		\rho &= \nabla\cdot R,\quad\text{where} \quad R \in L^2(Q_T)\quad \text{and} \quad R_t \in W_2^{0, \frac{l}{2}}(Q_T), \\
		\nabla\cdot u_0 &= \rho|_{t=0},\\
		d|_{t=0} &= 2[\mathbf{D}(v)\mathbf{n} - (\mathbf{D}(v)\mathbf{n}\cdot \mathbf{n})\mathbf{n}]|_{S_F},  \\
		d\cdot \mathbf{n} &= 0.
		\end{split}
		\end{equation*}
		Then system (\ref{stokes system}) with general initial data $v_0$ has a solution 
		\[
		v \in W_2^{l+2, \frac{l}{2}+1}(Q_T) \cap C_T W_2^{l+1}(\Omega), \quad q \in W_2^{l, \frac{l}{2}}(Q_T), \quad \nabla q \in W_2^{l, \frac{l}{2}}(Q_T),\quad q \in W_2^{l+ \frac{1}{2}, \frac{l}{2}+ \frac{1}{4}}(S_{F,T}),
		\]
		with the estimate
		\begin{equation}
		\begin{split}
		& \|v\|_{H^{l+2, \frac{l}{2}+1}(Q_T)} + \|q\|^{(l)}_{Q_T} + \|\nabla q\|^{(l)}_{Q_T} + \|q\|_{W_2^{l+ \frac{1}{2}, \frac{l}{2}+ \frac{1}{4}}(S_{F,T})}  \\
		&\leq C(T) \{ \|f\|^{(l)}_{Q_T} + \|\rho\|_{W_2^{l+1, \frac{l+1}{2}}(S_{F,T})} + \|R\|_{W_2^{0, \frac{l}{2}+1}(Q_T)} + T^{- \frac{l}{2}}\|R_t\|_{L^2(Q_T)}  \\
		&\quad + \|(b,d)\|_{W_2^{l+ \frac{1}{2}, \frac{l}{2}+ \frac{1}{4}}(S_{F,T})} + T^{- \frac{l}{2}}\|b\|_{W_2^{ \frac{l}{2},0}(S_{F,T})} + \|u_0\|_{W_2^{l+1}(\Omega)} \},
		\end{split}
		\end{equation}
		where $C(T)$ is time dependent constant on $T$ non-decreasingly.
	\end{proposition}
	
	\begin{proof}
		Let us write $w = v - U$, where $U$ satisfies (\ref{U ineq}), $U(0)=v_0$, and $U \in W^{l+2,\frac{l}{2}+1}(Q_{\infty})$. Then $(w,q)$ solves linear Stokes problem with zero initial data,
		\begin{equation}\label{w eq}
		\begin{cases}
		\frac{\partial w}{\partial t} -  \triangle w + \nabla q = f - \frac{\partial U}{\partial t} +  \triangle U := f^{\prime},\quad\text{in}\quad Q_T, \\
		\nabla\cdot w = \rho - \nabla\cdot U := \rho^{\prime},\quad \text{in}\quad Q_T, \\
		2 [\mathbf{D}(w)\mathbf{n} - (\mathbf{D}(w)\mathbf{n}\cdot\mathbf{n})\mathbf{n}] = d - 2 [\mathbf{D}(U)\mathbf{n} - (\mathbf{D}(U)\mathbf{n}\cdot\mathbf{n})\mathbf{n}] := d^{\prime},\quad \text{on}\quad S_{F}(t), \\
		q - 2 \mathbf{D}(w)\mathbf{n}\cdot\mathbf{n} := - b^{\prime},\quad \text{on}\quad S_{F}(t), \\
		w(t) = 0,\quad\text{on}\quad S_{B}, \\
		w(0) = 0,\quad\Omega\times\{t=0\}. \\
		\end{cases}
		\end{equation}
		We estimate $(f',\rho',d',b')$. Using definitions of $(f',\rho',d',b')$, 
		\begin{equation} \label{f' rho'}
		\begin{split}
		\|f'\|^{(l)}_{Q_T} &\leq C( \|f\|^{(l)}_{Q_T} + \|u_0\|_{W_2^{l+1}(\Omega)} ),  \\
		\|\rho'\|_{W_2^{l+1,0}(Q_T)} &\leq C ( \|\rho\|_{W_2^{l+1,0}(Q_T)} + \|u_0\|_{W_2^{l+1}(\Omega)} ).
		\end{split}
		\end{equation}
		At time $t=0$, we have $\nabla\cdot w_0 = 0 = \rho|_{t=0} - \nabla\cdot u_0 $ and $\rho^{\prime}|_{t=0} = 0$. If we write $\rho = \nabla\cdot R$, then
		$$
		\nabla\cdot w = \rho^{\prime} = \nabla\cdot(R-U) := \nabla\cdot R^{\prime},
		$$
		which gives,
		$$
		\|R^{\prime}\|_{W_2^{0, \frac{l}{2}+1}(Q_T)} \leq \|R\|_{W_2^{0, \frac{l}{2}+1}(Q_T)} + \|U\|_{W_2^{l+2, \frac{l}{2}+1}(Q_T)}.
		$$
		From Proposition~\ref{stokes with zero initial}, we should estimate, $\rho^{\prime}\equiv\nabla\cdot R^{\prime}$. It is easy to estimate $R^{\prime\prime}$ which satisfies $\rho^{\prime}=\nabla\cdot R^{\prime}=\nabla\cdot R^{\prime\prime}$ and $R^{\prime\prime}|_{t=0}=0$, because $H_{\gamma}^{l, \frac{l}{2}}$ requires compatibility condition. Let $\varphi$ solves Dirichlet Laplace problem, 
		$$
		\triangle\varphi = \nabla\cdot R' = \rho',\,\,\,\,\varphi|_{\partial\Omega} = 0,
		$$
		and we define $R^{\prime\prime} := \nabla\varphi$. Then 
		$$
		\rho^{\prime} = \nabla\cdot R^{\prime} = \nabla\cdot R^{\prime\prime} .
		$$
		Moreover, since initial condition is divergence-free, we can subtract initial condition from $R''$ while maintains divergence-free property. It gives $R''|_{t=0} = 0$. Now we use $R''$ instead of $R'$. We have an estimate for $R''$.
		\begin{equation} \label{R"}
		\begin{split}
		\|R''\|_{H_0^{0, \frac{l}{2}+1}(Q_T)} &= \|\nabla\varphi_t\|_{H_0^{0, \frac{l}{2}}(Q_T)} \leq C \|R'\|_{H_0^{0, \frac{l}{2}}(Q_T)}  \\
		&\leq C \left( \|R_t\|_{W_2^{0, \frac{l}{2}}(Q_T)} + T^{- \frac{l}{2}}\|R_t\|_{L^2(Q_T)} + \|u_0\|_{W_2^{l+1}(\Omega)} \right).
		\end{split}
		\end{equation}
		To estimate $d'$, we use condition $d'|_{t=0}=0$, Lemma \ref{solo lemma 63}, and the trace theorem to get
		\begin{equation} \label{d'}
		\begin{split}
		\|d'\|_{H_0^{l+ \frac{1}{2}, \frac{l}{2}+ \frac{1}{4}}(S_{F,T})} \leq C \left( \|d\|_{W_2^{l+ \frac{1}{2}, \frac{l}{2}+ \frac{1}{4}}(S_{F,T})} + \|u_0\|_{W_2^{l+1}(\Omega)} \right).
		\end{split}
		\end{equation}
		Similarly,
		\begin{equation} \label{b'}
		\begin{split}
		\|b'\|_{H_0^{l+ \frac{1}{2}, \frac{1}{2}, \frac{l}{2}}(S_{F,T})} \leq C \left( \|b\|_{W_2^{l+ \frac{1}{2}, \frac{l}{2}+ \frac{1}{4}}(S_{F,T})} + T^{- \frac{l}{2}}\|b\|_{W_2^{ \frac{1}{2},0}(S_{F,T})} + \|u_0\|_{W_2^{l+1}(\Omega)} \right).
		\end{split}
		\end{equation}
		
		\noindent We apply the result of Property \ref{stokes with zero initial} with (\ref{f' rho'}), (\ref{R"}), (\ref{d'}), and (\ref{b'}) to get, 
		\begin{equation} \label{w}
		\begin{split}
		&\|w\|_{H_0^{l+2, \frac{l}{2}+1}(Q_T)} + \|q\|_{H_0^{l+1,1, \frac{l}{2}}(Q_T)} \leq e^{\gamma T} \left( \|w\|_{H_{\gamma}^{l+2, \frac{l}{2}+1}(Q_T)} + \|\nabla q\|_{H_{\gamma}^{l, \frac{l}{2}}(Q_T)} \right)  \\
		&\leq C e^{\gamma T} \{ \|f'\|_{H_{\gamma}^{l, \frac{l}{2}}(Q_T)} + \|\rho'\|_{H_{\gamma}^{l+1, \frac{l+1}{2}}(Q_T)} + \|R''\|_{H_{\gamma}^{0, \frac{l}{2}+1}(Q_T)} + \|d'\|_{H_{\gamma}^{l+ \frac{1}{2}, \frac{l}{2}+ \frac{1}{4}}(S_{F,T})} + \|b'\|_{H_{\gamma}^{l+ \frac{1}{2}, \frac{1}{2}, \frac{l}{2}}(S_{F,T})} \}\\
		&\leq C(T) \{ \|f\|^{(l)}_{Q_T} + \|\rho\|_{W_2^{l+1, \frac{l+1}{2}}(Q_T)} + T^{- \frac{l}{2}}\|R_t\|_{L^2(Q_T)} + \|R\|_{W_2^{0, \frac{l}{2}+1}(Q_T)}  \\
		&\quad + \|(d,b)\|_{W_2^{l+ \frac{1}{2}, \frac{l}{2}+ \frac{1}{4}}(S_{F,T})} + T^{- \frac{l}{2}}\|b\|_{W_0^{ \frac{1}{2},0}(S_{F,T})} + \|v_0\|_{W_2^{l+2}(\Omega)} \},  \\
		\end{split}
		\end{equation}
		where $\gamma$ is a fixed constant which is found in Proposition~\ref{stokes with zero initial} and $C(T)$ depends on time $T$ non-decreasingly on $T$ which means it does not blow up as $T\rightarrow 0$. We also note that norms of $H^{r,\frac{r}{2}}_{\gamma}$ and $H^{r,\frac{r}{2}}_{0}$ are equivalent when $T < \infty$. Therefore (\ref{w eq}) yields above estimate for $\|w\|_{H_{\gamma}^{l+2, \frac{l}{2}+1}(Q_T)} + \|q\|_{H_{\gamma}^{l+1,1, \frac{l}{2}}(Q_T)}$. Meanwhile, from boundary condition,
		\begin{equation} \label{q}
		\|q\|_{W_2^{0,\frac{l}{2}+\frac{1}{4}}(S_{F,T})} \leq \|2 \mathbf{D}(v)\mathbf{n}\cdot \mathbf{n}\|_{W_2^{0,\frac{l}{2}+\frac{1}{4}}(S_{F,T})} + \|b\|_{W_2^{0,\frac{l}{2}+\frac{1}{4}}(S_{F,T})}.  \\
		\end{equation}
		Now we estimate $L_{T}^{\infty}$ type estimates for $w$ and $U$. From Proposition 4 in \cite{Lions} or Lemma 2.3 in \cite{JB}, we have interpolation
		\begin{equation}
			W^{r,\frac{r}{2}}(Q_{T}) \rightarrow H^{p}(0,T; W^{r-2p, \frac{r}{2}-p}).
		\end{equation}
		Since $w$ has zero initial data, we use Lemma \ref{extention lemma} to obtain
		\begin{equation} \label{wC}
		\begin{split}
		\|w\|_{C_{T} W_{2}^{l+1}(\Omega)} &\leq \|w_{ext}\|_{C_\infty W_{2}^{l+1}(\Omega)}
		\leq C \left( \|w_{ext}\|_{L^2_\infty W_2^{l+2}(\Omega)} + \|\partial_t w_{ext}\|_{L^2_\infty W_2^{l}(\Omega)} \right)  \\
		&\leq C \|w_{ext}\|_{W_2^{l+2,\frac{l}{2}+1}(Q_\infty)} \leq C \|w\|_{W_2^{l+1,\frac{l}{2}+1}(Q_T)},
		\end{split}
		\end{equation}
		where $w_{\text{ext}} \in W_2^{l+2, \frac{l}{2}+1}$ is an extension of $w$ which satisfies (we can choose time-independent $C$)
		$$
		\|w_{\text{ext}}\|_{W_2^{l+2,\frac{l}{2}+1}(Q_\infty)} \leq C \|w\|_{W_2^{l+2,\frac{l}{2}+1}(Q_T)}.
		$$ 
		Similarly,
		\begin{equation} \label{U}
		\begin{split}
		\|U\|_{C_T W_2^{l+1}(\Omega)} &\leq \|U_{\text{ext}}\|_{C_\infty W_2^{l+1}(\Omega)} \leq C \left( \|U_{\text{ext}}\|_{L^2_\infty W_2^{l+2}(\Omega)} + \|\partial_t U_{\text{ext}}\|_{L^2_\infty W_2^{l}(\Omega)} \right) \\
		&\leq C \|U_{\text{ext}}\|_{W_2^{l+2,\frac{l}{2}+1}(Q_\infty)} \leq C \|v_0\|_{W_2^{l+1}(\O)},
		\end{split}
		\end{equation}
		where $U_{\text{ext}} \in W_2^{l+2,\frac{l}{2}+1}$ is an extension of $U$ which satisfies
		\begin{equation*} \label{U ext}
		\|U_{\text{ext}}\|_{W_2^{l+2,\frac{l}{2}+1}(Q_\infty)} \leq C \|v_0\|_{W_2^{l+1}(\Omega)},
		\end{equation*}
		from Lemma \ref{extention lemma}. From (\ref{infty space}),
		$$
		(\|v\|^{(l+2)}_{Q_T})^{2} + \sup_{t<T} \|v\|^2_{W_2^{l+1}(\Omega)} \leq C\left\{ \|U\|^{(l+2)2}_{Q_T} + \sup_{t<T}\|U\|^2_{W_2^{l+1}(\Omega)} + \|w\|^{(l+2)2}_{Q_T} + \sup_{t<T}\|w\|^2_{W_2^{l+1}(\Omega)} \right\},
		$$
		so putting Lemma \ref{extention lemma}, (\ref{w}), (\ref{q}), (\ref{wC}), and (\ref{U}) together, we finish the proof.
	\end{proof}

	\section{constant coefficient nonlinear problem}
	In this section, we solve constant coefficient nonlinear problem,
	\begin{equation} \label{cst coeff nonlin}
	\begin{cases}
	v_t -  \triangle v + \nabla q = (B\cdot\nabla)B + f ,\quad\text{in}\quad Q_{T}, \\
	\nabla\cdot v = \rho,\quad \text{in}\quad Q_{T}, \\
	v(0) = v_0,\quad\text{in}\quad\Omega\times\{t=0\}, \\
	v = 0,\quad\text{on}\quad S_{B},  \\
	2 [\mathbf{D}(v)\mathbf{n} - (\mathbf{D}(v)\mathbf{n}\cdot \mathbf{n})\mathbf{n}] = d,\quad \text{on}\quad S_{F,T}, \\
	-q + 2 \mathbf{D}(v) \mathbf{n}\cdot \mathbf{n} = b \,\,\,\,\text{on}\,\,\,\,S_{F,T}, \\
	B_t -  \triangle B = (B\cdot\nabla)v + g,\quad\text{in}\quad Q_T ,\\
	B = 0,\quad\text{on}\quad  S_{F}\cup \O_{V} , \\
	B(0) = B_0 ,\quad \text{in}\quad \Omega\times\{t=0\}. \\
	B = 0,\quad\text{on}\quad S_{B}.  \\
	\end{cases}
	\end{equation}
	
	To control nonlinear terms $(B\cdot\nabla)B$ and $(B\cdot\nabla)v$, we prove the following Lemma.
	\begin{lemma} \label{nonlinear est}
		When $l\in(\frac{1}{2},1)$, we have the following nonlinear estimate.
		$$
		\|F\nabla G\|^{(l)}_{Q_T} \leq C ( T + T^{ \frac{1-l}{2}} + T^{ \frac{1}{2}} ) \|F\|_{H^{l+2, \frac{l}{2}+1}(Q_T)} \|G\|_{H^{l+2, \frac{l}{2}+1}(Q_T)}.
		$$
	\end{lemma}
	\begin{proof}
		\begin{equation*}
		\begin{split}
		\|F\nabla G\|^{(l)2}_{Q_T} &= \|F\nabla G\|^2_{W_2^{l, \frac{l}{2}}(Q_T)} + T^{-l}\|F\nabla G\|^2_{L^2(Q_T)}  \\
		&= \|F\nabla G\|^2_{W_2^{l,0}(Q_T)} + \|F\nabla G\|^2_{W_2^{0, \frac{l}{2}}(Q_T)} + T^{-l}\|F\nabla G\|^2_{L^2(Q_T)}.
		\end{split}
		\end{equation*}
		Using algebraic property of Lemma \ref{interpolation lemma} with $l+1>\frac{3}{2}$, 
		\begin{equation} \label{F delta G}
		\begin{split}
		\|F \nabla G\|^{(l)}_{Q_T} &\leq \|F \nabla G\|_{W_2^{l,0}(Q_T)} + \|F \nabla G\|_{W_2^{0, \frac{l}{2}}(Q_T)} + T^{- \frac{l}{2}}\|F \nabla G\|_{L^2(Q_T)}  \\
		&\leq C \Big\{ T \|F\|_{C_T W_2^{l+1}(\Omega)}\|G\|_{C_T W_2^{l+1}(\Omega)} + \underbrace{ T^{- \frac{l}{2}}\|F\nabla G\|_{L^2(Q_T)} }_{(II)} \Big\}  \\
		&\quad + \Big( \underbrace{ \int_0^T \int_0^T \frac{\left\| (F \nabla G)(t)-(F \nabla G)(s) \right\|^2_{L^2(\Omega)}}{\left| t-s \right|^{1+l}} dtds }_{(I)} \Big)^{ \frac{1}{2}} .
		\end{split}
		\end{equation}
		We focus on the last term $(I)$. Let us write $t-s = h$. Domain can be divided symmetrically into two regions $t>s$ and $s>t$. Then we change the order of integrals and apply Lemma~\ref{interpolation lemma},
		\begin{equation} \label{(I)}
		\begin{split}
		(I)&= \int_0^T \int_0^T \frac{\left\| (F \nabla G)(t)-(F \nabla G)(s) \right\|^2_{L^2(\Omega)}}{\left| t-s \right|^{1+l}} dtds  \\
		&\leq C \int_0^T \frac{dh}{h^{1+l}} \int_h^T \left\|F(s)\nabla (G(s+h)-G(s)) \right\|_{L^2(\Omega)}^2 ds   \\
		&\quad + C \int_0^T \frac{dh}{h^{1+l}} \int_h^T \left\|(F(s+h)-F(s))\nabla G(s)\right\|_{L^2(\Omega)}^2 ds \\
		&\leq C \|F\|^2_{C_T W_2^{l+1}(\Omega)}\int_0^T \frac{dh}{h^{1+l}} \int_h^T \left\|G(s+h) - G(s)) \right\|_{W_2^1(\Omega)}^2 ds  \\
		&\quad + C \|G\|_{C_T W_2^{l+1}(\Omega)}\int_0^T \frac{dh}{h^{1+l}} \int_h^T \left\|F(s+h)-F(s)\right\|_{W_2^1(\Omega)}^2 ds  \\
		&\leq CT \|F\|^2_{C_T W_2^{l+1}(\Omega)}\int_0^T \frac{dh}{h^{1+(l+1)}} \int_h^T \left\|G(s+h) - G(s)) \right\|_{W_2^1(\Omega)}^2 ds  \\
		&\quad + CT \|G\|_{C_T W_2^{l+1}(\Omega)}\int_0^T \frac{dh}{h^{1+(l+1)}} \int_h^T \left\|F(s+h)-F(s)\right\|_{W_2^1(\Omega)}^2 ds  \\
		&\leq CT \|F\|^2_{C_T W_2^{l+1}(\Omega)} \|G\|^2_{H^{1+2(\frac{l+1}{2}),\frac{1}{2}+\frac{l+1}{2}}(Q_T)} + CT \|G\|^2_{C_T W_2^{l+1}(\Omega)} \|F\|^2_{H^{1+2(\frac{l+1}{2}),\frac{1}{2}+\frac{l+1}{2}}(Q_T)}  \\
		&\leq CT \|F\|^2_{H^{l+2, \frac{l}{2}+1}(Q_T)} \|G\|^2_{H^{l+2, \frac{l}{2}+1}(Q_T)}.
		\end{split}
		\end{equation}
		For $(II)$, from Lemma~\ref{interpolation lemma} and $l\in(\frac{1}{2},1)$,
		\begin{equation} \label{(II)}
		\begin{split}
		(II) &= T^{- \frac{l}{2}}\|F\nabla G\|_{L^2(Q_T)} \\
		&\leq CT^{- \frac{l}{2}} \left( \int_0^T \|F\|^2_{W_2^{\frac{3}{2}-l}(\Omega)}\|\nabla G\|^2_{W_2^{l}(\Omega)} dt \right)^{ \frac{1}{2}}  \\
		&\leq C T^{ \frac{1-l}{2}} \|F\|_{C_T W_2^{l+1}(\Omega)} \|G\|_{C_T W_2^{l+1}(\Omega)}.
		\end{split}
		\end{equation}
		Combining (\ref{F delta G}), (\ref{(I)}), and (\ref{(II)}), we have the following estimate.
		\begin{equation*}
		\begin{split}
		\|F \nabla G\|^{(l)}_{Q_T} &\leq C \big\{ T \|F\|_{C_T W_2^{l+1}(\Omega)}\|G\|_{C_T W_2^{l+1}(\Omega)} + T^{ \frac{1-l}{2}} \|F\|_{C_T W_2^{l+1}(\Omega)} \|G\|_{C_T W_2^{l+1}(\Omega)}  \\
		&\quad + T^{ \frac{1}{2}} \|F\|_{H^{l+2, \frac{l}{2}+1}(Q_T)} \|G\|_{H^{l+2, \frac{l}{2}+1}(Q_T)} \big\}  \\
		&\leq C ( T + T^{ \frac{1-l}{2}} + T^{ \frac{1}{2}} ) \|F\|_{H^{l+2, \frac{l}{2}+1}(Q_T)} \|G\|_{H^{l+2, \frac{l}{2}+1}(Q_T)}.
		\end{split}
		\end{equation*}
	\end{proof}

	We use Proposition~\ref{prop general stokes} in section 3 and Lemma~\ref{nonlinear est} to solve system (\ref{cst coeff nonlin}). 
	\begin{proposition} \label{prop cst coeff nonlin}
		Let $l\in( \frac{1}{2},1)$, and $S_F \in W_2^{l+\frac{3}{2}}$. Assume that $(f,\rho,v_0,(b,d))$ in $(\ref{cst coeff nonlin})$ satisfy
		\begin{equation*}
		\begin{split}
		(f,\rho,u_0,(b,d)) &\in W_2^{l, \frac{l}{2}}(Q_T) \times W_2^{l+1, \frac{l+1}{2}}(Q_T) \times W_2^{l+1}(\Omega) \times W_2^{l+ \frac{1}{2}, \frac{l}{2}+ \frac{1}{4}}(S_{F,T}),  \\
		\rho &= \nabla\cdot R, \quad R \in L^2(Q_T), \quad 	R_t \in W_2^{0, \frac{l}{2}}(Q_T).
		\end{split}
		\end{equation*}
		Also assume that the following compatibility conditions hold,
		$$
		\nabla\cdot v_0 = \rho|_{t=0},\quad d|_{t=0} = 2 [\mathbf{D}(v)\mathbf{n} - (\mathbf{D}(v)\mathbf{n}\cdot \mathbf{n})\mathbf{n}]|_{S_F},\quad d\cdot \mathbf{n} = 0.
		$$
		Then system (\ref{cst coeff nonlin}) has a solution 
		\begin{equation*}
		\begin{split}
		v &\in W_2^{l+2, \frac{l}{2}+1}(Q_T) \ \cap \ C_T W_2^{l+1}(\Omega),\quad q \in W_2^{l, \frac{l}{2}}(Q_T), \\
		\nabla q &\in W_2^{l, \frac{l}{2}}(Q_T), \quad q \in W_2^{l+ \frac{1}{2}, \frac{l}{2}+ \frac{1}{4}}(S_{F,T}).
		\end{split}
		\end{equation*}
	\end{proposition}
	\begin{proof}
		We construct the following iteration scheme:
		\begin{equation}
		\begin{cases}
		v^{(m+1)}_t -  \triangle v^{(m+1)} + \nabla q^{(m+1)} = (B^{(m)}\cdot\nabla)B^{(m)} + f ,\quad\text{in}\quad Q_T, \\
		\nabla\cdot v^{(m+1)} = \rho,\quad\text{in}\quad Q_T, \\
		v^{(m+1)}(0) = v_0,\quad \text{in}\quad \Omega\times\{t=0\}, \\
		v^{(m+1)} = 0,\quad \text{on}\quad S_{B}, \\
		2 [\mathbf{D}(v^{(m+1)})\mathbf{n} - (\mathbf{D}(v^{(m+1)})\mathbf{n}\cdot \mathbf{n})\mathbf{n}] = d,\quad\text{on}\quad S_{F,T}, \\
		-q + 2 \mathbf{D}(v^{(m+1)})\mathbf{n}\cdot \mathbf{n} = b,\quad\text{on}\quad S_{F,T}, \\
		B_t^{(m+1)} -  \triangle B^{(m+1)} = (B^{(m)}\cdot\nabla)v^{(m)} + g ,\quad\text{in}\quad Q_T, \\
		B^{(m+1)} = 0, \quad\text{on}\quad S_{F,T}\cup \{\mathbb{R}^{3}\backslash Q_{T}\}, \\
		B^{(m+1)}(0) = B_0,\quad\text{in}\quad\Omega\times\{t=0\}, \\
		B^{(m+1)} = 0,\quad \text{on}\quad S_{B}, \\
		(v^{(0)},q^{(0)},B^{(0)}) = (0,0,0).
		\end{cases}
		\end{equation}
		
		For PDE for $B$, it is easy to estimate Heat equation with zero boundary condition. $\O$ is unbounded domain, but $\O$ is bounded in vertical direction, so we can use Poincar\'e inequality $\|u\|_{L^{2}(\O)} \leq \|\nabla u\|_{L^{2}(\O)}$. Therefore, from standard parabolic estimate with integer $l \in\mathbb{N}$,
		\[
		B_t - \triangle B = g,\quad\text{in}\quad Q_T, \\
		\]
		with Dirichlet boundary data has standard estimates,
		\begin{equation*}
		\begin{split}
		\|B\|^2_{L^{2}W_{2}^{l+2}} + \|B\|^2_{W_{2}^{l}L^{2}} + \|B\|_{C_{T}W_{2}^{l+1}} &\lesssim \|B_{0}\|_{W^{l+1}_{2}}^{2} + \|g\|^{2}_{W_{2}^{l}(Q_{T})}. 
		\end{split} 
		\end{equation*}	
		And also by interpolation, this holds for non-integer $l$ also. Since,
		\begin{equation*}
		\begin{split}
		\|B\|^2_{H^{l+2, \frac{l}{2}+1}(Q_T)} &:= \sup_{t<T}\|B(t)\|^2_{W_2^{l+1}(\Omega)} + \big( \|B\|^2_{W_2^{l, \frac{l}{2}}(Q_T)} +  \sum_{|s|=2} \|D_{x}^{s} B\|^2_{W_2^{l, \frac{l}{2}}(Q_T)} \big)  \\
		&\quad + T^{-l} \big( \|B\|^2_{L^2(Q_T)} + \sum_{|s|=2} \|D_{x}^{s} B\|^2_{L^2(Q_T)} \big) + \sum_{|s|=0}^{1} \|D_{x}^{s} B\|_{L^{2}(Q_{T})}^{2}, \\ 
		\end{split} 
		\end{equation*}	
		and $l\in(\frac{1}{2},1)$,
		\begin{equation} \label{Heat est}
		\begin{split}
		\|B\|^2_{H^{l+2, \frac{l}{2}+1}(Q_T)} &\leq T^{1-l} \big(  \|B\|^2_{L^{2}W_{2}^{l+2}} + \|B\|^2_{W_{2}^{l}L^{2}} + \|B\|_{C_{T}W_{2}^{l+1}}  \big)  \\
		&\lesssim C(T) \Big( \|B_{0}\|_{W^{l+1}_{2}}^{2} +  (\|g\|^{(l)}_{Q_T})^{2} \Big),  \\
		\end{split}
		\end{equation}
		where $C(T)$ depends on $T$ increasingly.  \\
		
		From Proposition~\ref{prop general stokes} and (\ref{Heat est}), we have a unique solution $(v^{(m+1)},q^{(m+1)},B^{(m+1)})$, for given data $(v^{(m)},B^{(m)})$,
		\begin{equation}
		\begin{split}
		v^{(m+1)} &\in W_2^{l+2, \frac{l}{2}+1}(Q_T)\cap C_T W_2^{l+1}(\Omega), \\
		B^{(m+1)} &\in W_2^{l+2, \frac{l}{2}+1}(Q_T)\cap C_T W_2^{l+1}(\Omega), \\
		q^{(m+1)} &\in W_2^{l, \frac{l}{2}}(Q_T),  \\
		\nabla q^{(m+1)} &\in W_2^{l, \frac{l}{2}}(Q_T),  \\
		q &\in W_2^{l, \frac{l}{2}}(S_{F,T}).
		\end{split}
		\end{equation}
		To get uniform bounds, we first define 
		\begin{equation} \label{mathcal A}
		\begin{split}
		\mathcal{A}^{(m+1)} &:= \|v^{(m+1)}\|_{H^{l+2, \frac{l}{2}+1}(Q_T)} + \|B^{(m+1)}\|_{H^{l+2, \frac{l}{2}+1}(Q_T)} + \|q^{(m+1)}\|^{(l)}_{Q_T} \\
		&\quad + \|\nabla q^{(m+1)}\|^{(l)}_{Q_T} + \|q^{(m+1)}\|_{W_2^{l+ \frac{1}{2}, \frac{l}{2}+ \frac{1}{4}}(S_{F,T})} .
		\end{split}
		\end{equation}

		Now, we use estimates of Proposition~\ref{prop general stokes} and (\ref{Heat est}) to get,
		\begin{equation} \label{A est}
		\begin{split}
		\mathcal{A}^{(m+1)} &\leq C(T) \{ \|(B^{(m)}\cdot\nabla)B^{(m)}\|^{(l)}_{Q_T} + \|(B^{(m)}\cdot\nabla)v^{(m)}\|^{(l)}_{Q_T} + \|f\|^{(l)}_{Q_T} + \|g\|^{(l)}_{Q_T}  \\
		&\quad + \|\rho\|_{W_2^{l+1, \frac{l+1}{2}}(S_{F,T})} + \|R\|_{W_2^{0, \frac{l}{2}+1}(Q_T)} + T^{- \frac{l}{2}}\|R\|_{L^2(Q_T)}  \\
		&\quad + \|(b,d)\|_{W_2^{l+ \frac{1}{2}, \frac{l}{2}+ \frac{1}{4}}(S_{F,T})} + T^{- \frac{l}{2}}\|b\|_{W_2^{ \frac{l}{2},0}(S_{F,T})} + \|u_0\|_{W_2^{l+1}(\Omega)} + \|H_0\|_{W_2^{l+1}(\Omega)}\}.
		\end{split}
		\end{equation}
		
		We also define data part as
		\begin{equation} \label{mathcal D}
		\begin{split}
		\mathcal{D}(T) &:= \|f\|^{(l)}_{Q_T} + \|g\|^{(l)}_{Q_T} + \|\rho\|_{W_2^{l+1, \frac{l+1}{2}}(S_{F,T})} + \|R\|_{W_2^{0, \frac{l}{2}+1}(Q_T)}  \\
		&\quad + T^{- \frac{l}{2}}\|R\|_{L^2(Q_T)} + \|(b,d)\|_{W_2^{l+ \frac{1}{2}, \frac{l}{2}+ \frac{1}{4}}(S_{F,T})} \\
		&\quad + T^{- \frac{l}{2}}\|b\|_{W_2^{ \frac{l}{2},0}(S_{F,T})} + \|v_0\|_{W_2^{l+1}(\Omega)} + \|B_0\|_{W_2^{l+1}(\Omega)}.
		\end{split}
		\end{equation}
		Using this definition, (\ref{A est}) is written by,
		\begin{equation} \label{mathcal A est}
		\begin{split}
		\mathcal{A}^{(m+1)} \leq  C(T)\left( \mathcal{D}(T) + ( T + T^{ \frac{1-l}{2}} + T^{ \frac{1}{2}} ) \mathcal{A}^{(m)2} \right).
		\end{split}
		\end{equation}
		Now we suffice to show uniform bound and contraction mapping to apply fixed point argument. \\
		$\mathbf{Uniform\,\,bound} \ $ By Cauchy sequence argument, we can pick sufficiently small $T_0 > 0$ such that
		\begin{equation} \label{unif bdd 1}
		\mathcal{A}^{(m)} \leq 2(C(0)+1)\mathcal{D}(T_0) \doteq M_{T_0},\,\,\,\,\forall m \in \mathbb{N}.
		\end{equation}
		$\mathbf{Contraction\,\,mapping} \ $ Let us define difference,
		$$
		v^{(m+1)} - v^{(m)} := \mathcal{V}^{(m+1)},\quad B^{(m+1)} - B^{(m)} := \mathcal{B}^{(m+1)},\quad q^{(m+1)} - q^{(m)} := \mathcal{Q}^{(m+1)}.
		$$
		Then $(\mathcal{V}^{(m+1)}, \mathcal{B}^{(m+1)}, \mathcal{Q}^{(m+1)})$ solves
		\begin{equation}
		\begin{cases}
		\mathcal{V}^{(m+1)}_t -  \triangle \mathcal{V}^{(m+1)} + \nabla \mathcal{Q}^{(m+1)} = (B^{(m)}\cdot\nabla)\mathcal{B}^{(m)} - (\mathcal{B}^{(m)}\cdot\nabla)B^{(m-1)},\quad\text{in}\quad Q_T, \\
		\nabla\cdot \mathcal{V}^{(m+1)} = 0,\quad \text{in}\quad Q_T, \\
		\mathcal{V}^{(m+1)}(0) = 0,\quad \text{in}\quad \Omega\times\{t=0\}, \\
		2 [\mathbf{D}(\mathcal{V}^{(m+1)})\mathbf{n} - (\mathbf{D}(\mathcal{V}^{(m+1)})\mathbf{n}\cdot \mathbf{n})\mathbf{n}] = 0,\quad \text{on}\quad S_{F,T}, \\
		-\mathcal{Q} + 2 \mathbf{D}(\mathcal{V}^{(m+1)})\mathbf{n}\cdot \mathbf{n} = 0 ,\quad \text{on}\quad S_{F,T}, \\
		\mathcal{B}_t^{(m+1)} -  \triangle \mathcal{B}^{(m+1)} = (B^{(m)}\cdot\nabla)\mathcal{V}^{(m)} - (\mathcal{B}^{(m)}\cdot\nabla)v^{(m-1)},\quad \text{in}\quad Q_T, \\
		\mathcal{B}^{(m+1)} = 0,\quad\text{on}\quad S_{F}\cup Q_{V}, \\
		\mathcal{V}^{(m+1)} = \mathcal{B}^{(m+1)} = 0,\quad\text{on}\quad S_{B}\\
		\mathcal{B}^{(m+1)}(0) = 0,\quad\text{in}\quad\Omega\times\{t=0\}. \\
		\end{cases}
		\end{equation}
		Again, using Proposition~\ref{prop general stokes}, (\ref{Heat est}), and Lemma~\ref{nonlinear est},
		\begin{equation}
		\begin{split}
		\bar{\mathcal{A}}^{(m+1)} &:= \|\mathcal{V}^{(m+1)}\|_{H^{l+2, \frac{l}{2}+1}(Q_T)} + \|\mathcal{B}^{(m+1)}\|_{H^{l+2, \frac{l}{2}+1}(Q_T)}  \\
		&\quad + \|\mathcal{Q}^{(m+1)}\|^{(l)}_{Q_T} + \|\nabla \mathcal{Q}^{(m+1)}\|^{(l)}_{Q_T} + \|\mathcal{Q}^{(m+1)}\|_{W_2^{l+ \frac{1}{2}, \frac{l}{2}+ \frac{1}{4}}(S_{F,T})}   \\
		&\leq C(T) \{ \|(B^{(m)}\cdot\nabla)\mathcal{B}^{(m)}\|^{(l)}_{Q_T} + \|(\mathcal{B}^{(m)}\cdot\nabla)B^{(m-1)}\|^{(l)}_{Q_T}  \\
		&\quad + \|(B^{(m)}\cdot\nabla)\mathcal{V}^{(m)}\|^{(l)}_{Q_T} + \|(\mathcal{B}^{(m)}\cdot\nabla)v^{(m-1)}\|^{(l)}_{Q_T} \}  \\
		&\leq C(T)( T + T^{ \frac{1-l}{2}} + T^{ \frac{1}{2}} ) M_{T_0} \bar{\mathcal{A}}^{(m)}.
		\end{split}
		\end{equation}
		We can find sufficiently small $T_1 > 0$ (Without loss of generality, we pick this so that smaller than $T_0$), such that
		$$
		C(t)( t + t^{ \frac{1-l}{2}} + t^{ \frac{1}{2}} ) M_{T_0} < 1,\,\,\,\,\forall t < T_1,
		$$
		since $C(T)$ depends on $T$ non-decreasingly, so that it is finite near $T=0$. Hence we have an unique solution via fixed point argument.
	\end{proof}
	
	\section{Proof of theorem \ref{main theorem}}
	In this section, we finish the proof of Theorem \ref{main theorem}. In subsection 5.1, we solve (\ref{Larg system}) and (\ref{Larg compatibility}). In subsection 5.2, propagation of divergence-free condition of magnetic field from initial data will be justified. 
	\subsection{Fully nonlinear system}
	In this subsection, we solve
	\begin{equation} \label{fully nonlin system}
	\begin{cases}
	v_t -  \triangle_v v + \nabla_v q = B\cdot\nabla_v B,\quad\text{in}\quad Q_T, \\
	B_t -  \triangle_v B = B\cdot\nabla_v v,\quad\text{in}\quad Q_T, \\
	\nabla_v\cdot v = 0,\quad\text{in}\quad Q_T, \\
	v(0) = v_0,\quad\Omega\times\{t=0\}, \\
	B(0) = B_0,\quad\Omega\times\{t=0\}, \\
	B = 0,\quad\text{on}\quad  S_{F,T}\cup \O_{V}, \\
	v = B = 0,\quad\text{on}\quad S_{B}, \\
	q - 2 \mathbf{D}_v(v)\mathbf{n}^{(v)}\cdot\mathbf{n}^{(v)} = gh,\quad\text{on}\quad S_{F,T}, \\
	2 \mathbf{D}(v)\mathbf{n}^{(v)} - 2 (\mathbf{D}_v(v)\mathbf{n}^{(v)}\cdot\mathbf{n}^{(v)})\mathbf{n}^{(v)} = 0,\quad\text{on}\quad S_{F,T}. \\
	\end{cases}
	\end{equation}
	
	Note that this system does not contain divergence-free condition for $B$, since (\ref{cst coeff nonlin}) does not include any condition about divergence-free for $B$. First we state a lemma which is similar as Lemma~\ref{nonlinear est}.  
	
	\begin{lemma} \label{fully nonlin est}
		For $l>\frac{1}{2}$, we have the following nonlinear estimate.
		$$
		\|FG\|^{(l)}_{Q_T} \leq C(T) \|F\|_{H^{l+1, \frac{l+1}{2}}(Q_T)} \|G\|_{H^{l+1, \frac{l+1}{2}}(Q_T)},
		$$
		where $C(T)$ depends on $T$ non-decreasingly.
	\end{lemma}
	\begin{proof}
		\begin{equation}
		\begin{split}
		\|FG\|^{(l)2}_{Q_T} &:= \|FG\|^2_{W_2^{l, \frac{l}{2}}(Q_T)} + T^{-l}\|FG\|^2_{L^2(Q_T)}  \\
		&= \|FG\|^2_{W_2^{l,0}(Q_T)} + \|FG\|^2_{W_2^{0, \frac{l}{2}}(Q_T)} + T^{-l}\|FG\|^2_{L^2(Q_T)},
		\end{split}
		\end{equation}
		\begin{equation}
		\begin{split}
		\|FG\|^{(l)}_{Q_T} &\leq \|FG\|_{W_2^{l,0}(Q_T)} + \|FG\|_{W_2^{0, \frac{l}{2}}(Q_T)} + T^{- \frac{l}{2}}\|FG\|_{L^2(Q_T)}  \\
		&\leq C \left\{ \|F\|_{L^2_T W_2^{l+1}(\Omega)}\|G\|_{C_T W_2^{l}(\Omega)} + T^{- \frac{l}{2}}\|F\|_{L^2_T W_2^1(\Omega)} \|G\|_{C_T W_2^l(\Omega)} \right\}  \\
		&+ \left( \underbrace{ \int_0^T \int_0^T \frac{\left\| (FG)(t)-(FG)(s) \right\|^2_{L^2(\Omega)}}{\left| t-s \right|^{1+l}} dtds }_{(III)} \right)^{ \frac{1}{2}},  \\
		\end{split}
		\end{equation}
		\begin{equation}
		\begin{split}
		(III) &:= \int_0^T \int_0^T \frac{\left\| (FG)(t)-(FG)(s) \right\|^2_{L^2(\Omega)}}{\left| t-s \right|^{1+l}} dtds  \\
		&\leq C \int_0^T \frac{dh}{h^{1+l}} \int_h^T \left\|F(s)(G(s+h)-G(s)) \right\|_{L^2(\Omega)}^2 ds + C \int_0^T \frac{dh}{h^{1+l}} \int_h^T \left\|(F(s+h)-F(s))G(s)\right\|_{L^2(\Omega)}^2 ds  \\
		&\leq C \left( \|F\|^2_{C_T W_2^{l}(\Omega)}\|G\|_{H^{l+1, \frac{l+1}{2}}(Q_T)} + \|G\|^2_{C_T W_2^{l}(\Omega)}\|F\|_{H^{l+1, \frac{l+1}{2}}(Q_T)} \right).
		\end{split}
		\end{equation}
		Using Lemma~\ref{interpolation lemma} ($l>\frac{1}{2}$), we finish the proof.
	\end{proof}
	
	Now, from above sections, we have a unique solution $(\mathbf{v},\mathbf{B},\mathbf{q})$ for $0 \leq t < T_1$ to system (\ref{cst coeff nonlin}) with $(f,\rho,u_0,d,b) = (0,0,0,0,-gh)$. We find a solution of the form 
	\[
	(v,B,q)=(\mathbf{v}+v^*,\mathbf{B}+B^*,\mathbf{q}+q^*).
	\]
	Then system (\ref{fully nonlin system}) becomes,
	\begin{equation} \label{perturb}
	\begin{cases}
	v^*_t -  \triangle_v v^* + \nabla_v q^* =  (\triangle_v - \triangle)\mathbf{v} - (\nabla_v - \nabla)\mathbf{q} + \mathbf{B}\cdot(\nabla_v - \nabla)\mathbf{B}, \\
	\,\,\,\,\,\,\,\,\,\,\,\,\,\,\,\,\,\,\,\,\,\,\,\,\,\,\,\,\,\,\,\,\,\,\,\,\,\,\,\,\,\,\,\,\,\,\,\,\,\,\,\,\,\,\,\,\,\,\,\, + \mathbf{B}\cdot\nabla_v B^* + B^*\cdot\nabla_v \mathbf{B} + B^*\cdot\nabla_v B^*,\quad\text{in}\quad Q_T, \\ 
	B^*_t -  \triangle_v B^* =  (\triangle_v - \triangle)\mathbf{B} + \mathbf{B}\cdot(\nabla_v - \nabla)\mathbf{v}, \\
	\,\,\,\,\,\,\,\,\,\,\,\,\,\,\,\,\,\,\,\,\,\,\,\,\,\,\,\,\,\,\,\,\,\,\,\,\,\,\,\,\,+ \mathbf{B}\cdot\nabla_v v^* + B^*\cdot\nabla_v \mathbf{v} + B^*\cdot\nabla_v v^*,\quad\text{in}\quad Q_T, \\ 
	\nabla_v\cdot v^* = -\nabla_v\cdot \mathbf{v},\quad\text{in}\quad Q_T, \\
	v^*(0) = v_0,\quad B^*(0) = B_0,\quad\text{in}\quad \Omega\times\{t=0\}, \\
	B^* = 0,\quad\text{on}\quad S_{F,T}\cup \O_{V}, \\
	v^* = B^* = 0,\quad\text{on}\quad S_{B}, \\
	-q^* + 2 \mathbf{D}_v(v^*)\mathbf{n}^{(v)}\cdot\mathbf{n}^{(v)} = -2 \mathbf{D}_v(\mathbf{v})\mathbf{n}^{(v)}\cdot\mathbf{n}^{(v)} + 2 \mathbf{D}_v(\mathbf{v}) \mathbf{n}\cdot \mathbf{n},\quad\text{on}\quad S_{F,T}, \\
	2 [\mathbf{D}(v^*)\mathbf{n}^{(v)} - (\mathbf{D}_v(v^*)\mathbf{n}^{(v)}\cdot\mathbf{n}^{(v)})\mathbf{n}^{(v)}] = -2 [\mathbf{D}(\mathbf{v})\mathbf{n}^{(v)} - (\mathbf{D}_v(\mathbf{v})\mathbf{n}^{(v)}\cdot\mathbf{n}^{(v)})\mathbf{n}^{(v)}],\quad\text{on}\quad S_{F,T}. \\
	\end{cases}
	\end{equation}
	
	We make the following iteration scheme to solve above (\ref{perturb}). Note that $\nabla\cdot \mathbf{v} = 0$.
	
	\begin{equation} \label{perturb iteration}
	\begin{cases}
	v^{*(m+1)}_t -  \triangle v^{*(m+1)} + \nabla q^{*(m+1)} =  (\triangle_m - \triangle)v^{(m)} - (\nabla - \nabla_m)q^{(m)} + \mathbf{B}\cdot(\nabla_m - \nabla)\mathbf{B}, \\
	\quad\quad\quad\quad\quad\quad\quad\quad\quad\quad\quad\quad\quad\quad\quad + \mathbf{B}\cdot\nabla_m B^{*(m)} + B^{*(m)}\cdot\nabla_m \mathbf{B} + B^{*(m)}\cdot\nabla_m B^{*(m)} \\
	\quad\quad\quad\quad\quad\quad\quad\quad\quad\quad\quad\quad\quad\quad\quad := f^{(m)},\quad\text{in}\quad Q_T, \\ 
	B^{*(m+1)}_t -  \triangle B^{*(m+1)} =  (\triangle_m - \triangle)B^{(m)} + \mathbf{B}\cdot(\nabla_m - \nabla)\mathbf{v}, \\
	\quad\quad\quad\quad\quad\quad\quad\quad\quad\quad + \mathbf{B}\cdot\nabla_m v^{*(m)} + B^{*(m)}\cdot\nabla_m \mathbf{v} + B^{*(m)}\cdot\nabla_m v^{*(m)} \\
	\quad\quad\quad\quad\quad\quad\quad\quad\quad\quad := g^{(m)},\quad\text{in}\quad Q_T, \\ 
	\nabla \cdot v^{*(m+1)} = (\nabla-\nabla_m)\cdot v^{(m)} := \rho^{(m)},\quad\text{in}\quad Q_T, \\
	v^{*(m+1)}(0) = u_0,\quad B^{*(m+1)}(0) = H_0,\quad\Omega\times\{t=0\}, \\
	B^{*(m+1)} = 0,\quad\text{on}\quad S_{F,T}\cup \O_{V}, \\
	v^{*(m+1)} = B^{*(m+1)} = 0,\quad\text{on}\quad S_{B}, \\
	-q^{*(m+1)} + 2 \mathbf{D}(v^{*(m+1)})\mathbf{n}\cdot \mathbf{n} = 2  [\mathbf{D}(v^{(m)})\mathbf{n}\cdot \mathbf{n} - \mathbf{D}_m(v^{(m)})\mathbf{n}^{(m)}\cdot \mathbf{n}^{(m)}]   \\
	\quad\quad\quad\quad\quad\quad\quad\quad\quad\quad\quad\quad\quad\quad := b^{(m)},\quad\text{on}\quad S_{F,T}, \\
	2 [\mathbf{D}(v^{*(m+1)})\mathbf{n} - (\mathbf{D}_v(v^{*(m+1)})\mathbf{n}\cdot \mathbf{n})\mathbf{n}] = -2  \{ [\mathbf{D}_m(v^{(m)})\mathbf{n}^{(m)} - (\mathbf{D}_m(v^{(m)})\mathbf{n}^{(m)}\cdot \mathbf{n}^{(m)})\mathbf{n}^{(m)}], \\
	\quad\quad\quad\quad\quad\quad\quad\quad\quad\quad\quad\quad\quad\quad\quad\quad\quad - [\mathbf{D}(v^{*(m)})\mathbf{n} - (\mathbf{D}(v^{*(m)})\mathbf{n} \cdot \mathbf{n})\mathbf{n}] \} \\
	\quad\quad\quad\quad\quad\quad\quad\quad\quad\quad\quad\quad\quad\quad\quad\quad\quad := d^{(m)},\quad\text{on}\quad S_{F,T} ,\\
	(v^{*(0)},B^{*(0)},q^{*(0)}) = (0,0,0),
	\end{cases}
	\end{equation}
	where
	\begin{equation}
	\begin{split}
	(v^{(m)},B^{(m)},q^{(m)}) &:= (\mathbf{v} + v^{*(m)},\mathbf{B} + B^{*(m)},\mathbf{q} + q^{*(m)}),  \\
	\nabla_m &:= \nabla_{v^{(m)}},\quad \mathbf{D}_m := \mathbf{D}_{v^{(m)}},\quad n^{(m)} := n^{v^{(m)}},  \\
	\rho^{(m)} &:= \nabla\cdot R^{(m)},\quad R^{(m)} := (\mathbf{I} - \mathcal{G}^{(m)})v^{(m)},\quad \mathcal{G}^{(m)} := \mathcal{G}^{v^{(m)}} ,
	\end{split}
	\end{equation}
	and 
	\begin{equation} \label{delta 12}
	T^{ \frac{1}{2}} \|v^{(m)}\|^{(l+2)}_{Q_{T}} \leq \delta_1,\,\,\,\,T^{ \frac{1}{2}} \|B^{(m)}\|^{(l+2)}_{Q_{T}} \leq \delta_2.
	\end{equation}
	Using Proposition~\ref{prop general stokes} and (\ref{Heat est}),
	\begin{equation}
	\begin{split}
	&\|v^{*(m+1)}\|_{H^{l+2, \frac{l}{2}+1}(Q_T)} + \|B^{*(m+1)}\|_{H^{l+2, \frac{l}{2}+1}(Q_T)}  \\
	&\quad + \|q^{*(m+1)}\|^{(l)}_{Q_T} + \|\nabla q^{*(m+1)}\|^{(l)}_{Q_T} + \|q^{*(m+1)}\|_{W_2^{l+ \frac{1}{2}, \frac{l}{2}+ \frac{1}{4}}(S_{F,T})}  \\
	&\leq C_*(T) \{ \|f^{(m)}\|^{(l)}_{Q_T} + \|g^{(m)}\|^{(l)}_{Q_T} + \|\rho^{(m)}\|_{W_2^{l+1, \frac{l+1}{2}}(S_{F,T})} + \|R^{(m)}\|_{W_2^{0, \frac{l}{2}+1}(Q_T)}  \\
	&\quad + T^{- \frac{l}{2}}\|R_t^{(m)}\|_{L^2(Q_T)} + \|(b^{(m)},d^{(m)})\|_{W_2^{l+ \frac{1}{2}, \frac{l}{2}+ \frac{1}{4}}(S_{F,T})}  \\
	&\quad + T^{- \frac{l}{2}}\|b^{(m)}\|_{W_2^{ \frac{l}{2},0}(S_{F,T})} + \|u_0\|_{W_2^{l+1}(\Omega)} + \|H_0\|_{W_2^{l+1}(\Omega)} \}.
	\end{split}
	\end{equation}
	
	Similar as section 4, we claim uniform bound and contraction mapping property.  \\
	$\mathbf{Uniform\,\,bound}$ 
	We can use Lemma~\ref{larg lemma}, with $u=u^{(m)}$ and $u'=0$. Especially, we compute laplacian as
	\begin{equation}
	\begin{split}
		& [\mathcal{G}^{(m)}]^{t} \nabla \big( \mathcal{G}^{(m)} \nabla v^{(m)} \big)  \\
		& = [\mathcal{G}^{(m)}]^{t} \big( \nabla \mathcal{G}^{(m)} \nabla v^{(m)} + \mathcal{G}^{(m)}\nabla^{2} v^{(m)} \big)  \\
		&= [\mathcal{G}^{(m)}]^{t} \mathcal{G}^{(m)} \nabla^{2} v^{(m)} + [\mathcal{G}^{(m)}]^{t} \nabla \mathcal{G}^{(m)} \nabla v^{(m)}  \\
		&= \Big( ([\mathcal{G}^{(m)}]^{t} - I)(\mathcal{G}^{(m)} - I) + ([\mathcal{G}^{(m)}]^{t} - \mathbf{I}) + (\mathcal{G}^{(m)} - \mathbf{I}) \Big) \nabla^{2} v^{(m)} + [\mathcal{G}^{(m)}]^{t} \nabla \mathcal{G}^{(m)} \nabla v^{(m)}.  \\
	\end{split}
	\end{equation}
	To treat second term, which includes the highest order term, we use Lemma~\ref{interpolation lemma} to control
	\begin{equation} \label{DG}
	\begin{split}
		\| [\mathcal{G}^{(m)}]^{t} \nabla \mathcal{G}^{(m)} \nabla v^{(m)} \|_{W^{l}_{2}} \lesssim \| [\mathcal{G}^{(m)}]^{t} \nabla \mathcal{G}^{(m)} \|_{W^{l}_{2}} \|\nabla v^{(m)} \|_{W^{l+1}_{2}},
	\end{split}
	\end{equation}
	since $l+1 > \frac{3}{2} = \frac{n}{2}$. Now we define 
	\begin{equation} \label{Z_m}
	\begin{split}
	Z_m &:= \|v^{(m)}\|_{H^{l+2, \frac{l}{2}+1}(Q_T)} + \|B^{(m)}\|_{H^{l+2, \frac{l}{2}+1}(Q_T)} + \|q^{(m)}\|^{(l)}_{Q_T} \\
	&+ \|\nabla q^{(m)}\|^{(l)}_{Q_T} + \|q^{(m)}\|_{W_2^{l+ \frac{1}{2}, \frac{l}{2}+ \frac{1}{4}}(S_{F,T})},
	\end{split}
	\end{equation}
	and use Lemma~\ref{nonlinear est}, (\ref{DG}) to estimate,    \\
	\begin{equation} \label{f_m}
	\begin{split}
	\|f^{(m)}\|^{(l)}_{Q_T} &\leq   \| \{ ([\mathcal{G}^{(m)}]^t-\mathbf{I})(\mathcal{G}^{(m)}-\mathbf{I}) + (\mathcal{G}^{(m)}-\mathbf{I})  \\
	&\quad + ([\mathcal{G}^{(m)}]^t-\mathbf{I}) \} \nabla\cdot\nabla v^{(m)} + \{ [\mathcal{G}^{(m)}\nabla]^t \cdot \mathcal{G}^{(m)} \} Dv^{(m)} \|^{(l)}_{Q_T} \\
	&\quad + \| (\mathcal{G}^{(m)}-\mathbf{I})\nabla q^{(m)}\|^{(l)}_{Q_T} + \| \mathbf{B}(\mathcal{G}^{(m)}-\mathbf{I})\nabla B^{*(m)} \|^{(l)}_{Q_T} \\
	&\quad + \| B^{*(m)}(\mathcal{G}^{(m)}-\mathbf{I})\nabla \mathbf{B} \|^{(l)}_{Q_T} + \|B^{*(m)}(\mathcal{G}^{(m)}-\mathbf{I})\nabla B^{*(m)} \|^{(l)}_{Q_T}  \\
	&\quad + \| \mathbf{B}\cdot\nabla B^{*(m)} \|^{(l)}_{Q_T} + \| B^{*(m)}\cdot\nabla \mathbf{B} \|^{(l)}_{Q_T} + \|  B^{*(m)}\cdot\nabla B^{*(m)} \|^{(l)}_{Q_T}  \\
	&\leq C(T,\delta_1,\delta_2) Z_m ( 1 + Z_m ) + \| \mathbf{B}\cdot\nabla B^{*(m)} \|^{(l)}_{Q_T} + \| B^{*(m)}\cdot\nabla \mathbf{B} \|^{(l)}_{Q_T} + \|  B^{*(m)}\cdot\nabla B^{*(m)} \|^{(l)}_{Q_T}  \\
	&\leq C(T,\delta_1,\delta_2) Z_m ( 1 + Z_m ) + C ( T + T^{ \frac{1-l}{2}} + T^{ \frac{l}{2}} ) Z_m^2  \\
	&\leq C(T,\delta_1,\delta_2) Z_m ( 1 + Z_m ),
	\end{split}
	\end{equation}
	and $C(T,\delta_1,\delta_2)$ is positive constant depending increasingly on its arguments satisfying $ C(T,\delta_1,\delta_2)\rightarrow 0$ as $T\rightarrow 0$.
	Using exactly same argument, we have the same estimate for $g^{(m)}$, $\rho^{(m)}$, $ R^{(m)}$, $R_t^{(m)}$, $d^{(m)}$, and $b^{(m)}$, 
	\begin{equation} \label{others}
	\begin{split}
	\|g^{(m)}\|^{(l)}_{Q_T} &\leq C(T,\delta_1,\delta_2) Z_m ( 1 + Z_m ),  \\
	\|\rho^{(m)}\|_{W_2^{l+1, \frac{l+1}{2}}(S_{F,T})} &\leq C(T,\delta_1,\delta_2) Z_m ( 1 + Z_m ),  \\
	\|R^{(m)}\|_{W_2^{0, \frac{l}{2}+1}(Q_T)} &\leq C(T,\delta_1,\delta_2) Z_m ( 1 + Z_m ) \\
	&\quad  + C(T,\delta_1,\delta_2) Z_m ( \|v^{(m)}\|_{W_2^{l+1,0}(Q_T)} + \|Dv^{(m)}\|_{W_2^{0, \frac{l}{2}}(Q_T)} ),  \\
	T^{- \frac{l}{2}}\|R_t^{(m)}\|_{L^2(Q_T)} &\leq C(T,\delta_1,\delta_2)T^{ \frac{1-l}{2}} Z_m ( 1 + Z_m ) \leq C(T,\delta_1,\delta_2) Z_m ( 1 + Z_m ),  \\
	\|d^{(m)}\|_{W_2^{l+ \frac{1}{2}, \frac{l}{2}+ \frac{1}{4}}(S_{F,T})} &\leq C(T,\delta_1,\delta_2) Z_m ( 1 + Z_m ), \\
	\|b^{(m)}\|_{W_2^{l+ \frac{1}{2}, \frac{l}{2}+ \frac{1}{4}}(S_{F,T})} & + T^{- \frac{l}{2}}\|b^{(m)}\|_{W_2^{ \frac{l}{2},0}(S_{F,T})} \leq C(T,\delta_1,\delta_2) Z_m ( 1 + Z_m ).
	\end{split}
	\end{equation}
	Similar as (\ref{Z_m}), we define,
	\begin{equation} \label{Z_m *}
	\begin{split}
	Z_m^* &:= \|v^{*(m)}\|_{H^{l+2, \frac{l}{2}+1}(Q_T)} + \|B^{*(m)}\|_{H^{l+2, \frac{l}{2}+1}(Q_T)} + \|q^{*(m)}\|^{(l)}_{Q_T} \\
	&+ \|\nabla q^{*(m)}\|^{(l)}_{Q_T} + \|q^{*(m)}\|_{W_2^{l+ \frac{1}{2}, \frac{l}{2}+ \frac{1}{4}}(S_{F,T})}.
	\end{split}
	\end{equation}
	Then, combining (\ref{f_m}) and (\ref{others}), 
	$$
	Z_{m+1}^* \leq h_0 + h_1 Z^*_m + h_2 Z_m^{*2},
	$$
	with positive constant $h_0,h_1,h_2$ with following properties.\\
	1) $h_0 = h_0(T,\delta_1,\delta_2)$ is monotone increasing function with all its argument.\\
	2) $h_{1,2} = h_{1,2}(T,\delta_1,\delta_2) \rightarrow 0$ as $T\rightarrow 0$.\\
	From Cauchy sequence argument, there exist $z^*$ such that if $Z_m^* \leq z^*$, then
	$$
	Z_{m+1}^* \leq h_0 + h_1 z^*_m + h_2 (z_m)^{*2} \leq z^*.
	$$
	Hence we have uniform bound,
	\begin{equation} \label{uniform bound 2}
	Z_m^* \leq z^*,\,\,\,\,\forall m
	\end{equation}
	$\mathbf{Contraction\,\,mapping}$ We use the similar notation in section 4 to denote, 
	$$
	v^{*(m+1)} - v^{*(m)} \doteq \mathcal{V}^{*(m+1)},\,\,\,\,B^{*(m+1)} - B^{*(m)} \doteq \mathcal{B}^{*(m+1)},\,\,\,\,q^{*(m+1)} - q^{*(m)} \doteq \mathcal{Q}^{*(m+1)}.
	$$
	For differential operators with for each $v^{m}$, we write
	\[
		\nabla_{m} := \nabla_{v^{m}},\quad \Delta_{m} := \Delta_{v^{m}}.
	\]
	Using (\ref{perturb iteration}) for $m+1$ and $m$, 
	\begin{equation}
	\begin{cases}
	\mathcal{V}^{*(m+1)}_t -  \triangle \mathcal{V}^{*(m+1)} + \nabla \mathcal{Q}^{*(m+1)} =   (\triangle_m - \triangle_{m-1})v^{(m)} +   (\triangle_{m-1} - \triangle )\mathcal{V}^{*(m)}, \\
	\,\,\,\,\,\,\,\,\,\,\,\,\,\,\,\,\,\,\,\,\,\,\,\,\,\,\,\,\,\,\,\,\,\,\,\,\,\,\,\,\,\,+ (\triangle_{m-1} - \triangle_m)q^{(m)} + (\triangle - \triangle_{m-1} )\mathcal{Q}^{*(m)} + \mathbf{B}\cdot(\nabla_m - \nabla_{m-1})\mathbf{B}, \\
	\,\,\,\,\,\,\,\,\,\,\,\,\,\,\,\,\,\,\,\,\,\,\,\,\,\,\,\,\,\,\,\,\,\,\,\,\,\,\,\,\,\,+ \mathbf{B}\cdot(\nabla_m - \nabla_{m-1})B^{*(m)} + \mathbf{B}\cdot\nabla_{m-1}\mathcal{B}^{*(m)}, \\
	\,\,\,\,\,\,\,\,\,\,\,\,\,\,\,\,\,\,\,\,\,\,\,\,\,\,\,\,\,\,\,\,\,\,\,\,\,\,\,\,\,\,+ B^{*(m)}\cdot(\nabla_m - \nabla_{m-1})\mathbf{B} + {\mathcal{B}}^{*(m)}\cdot\nabla_{m-1}\mathbf{B}, \\
	\,\,\,\,\,\,\,\,\,\,\,\,\,\,\,\,\,\,\,\,\,\,\,\,\,\,\,\,\,\,\,\,\,\,\,\,\,\,\,\,\,\, + B^{*(m)}\cdot\nabla_m \mathcal{B}^{*(m)} + B^{*(m)}\cdot( \nabla_m - \nabla_{m-1} )B^{*(m-1)} + \mathcal{B}^{*(m)}\cdot\nabla_{m-1} B^{*(m-1)} \\
	\,\,\,\,\,\,\,\,\,\,\,\,\,\,\,\,\,\,\,\,\,\,\,\,\,\,\,\,\,\,\,\,\,\,\,\,\,\,\,\,\,\, := f^{*(m)},\quad\text{in}\quad Q_T, \\
	\nabla\cdot\mathcal{V}^{*(m+1)} = (\nabla_{m-1}-\nabla_m)\cdot v^{(m)} + (\nabla - \nabla_{m-1})\cdot\mathcal{V}^{*(m)} := \rho^{*(m)}, \\
	\mathcal{V}^{*(m+1)}(0) = 0 \,\,\,\,\text{in}\,\,\,\,\Omega\times\{t=0\}, \\
	2 [\mathbf{D}(\mathcal{V}^{*(m+1)}) \mathbf{n} - (\mathbf{D}(\mathcal{V}^{*(m+1)})\mathbf{n}\cdot \mathbf{n})\mathbf{n}], \\
	\,\,\,\,\,\,\,\,\,\,\,\,\,\,\,\,\,\,\,\,\,\,\,\,\,\,\,\,\,\,\,\,\,\,\,\,\,\,\,\,\,\,= 2  \{[ \mathbf{D}(\mathcal{V}^{*(m)})\mathbf{n} - ( \mathbf{D}_m(v^{(m)})\mathbf{n}^{(m)} - \mathbf{D}_{m-1}(v^{(m-1)})\mathbf{n}^{(m-1)} )], \\
	\,\,\,\,\,\,\,\,\,\,\,\,\,\,\,\,\,\,\,\,\,\,\,\,\,\,\,\,\,\,\,\,\,\,\,\,\,\,\,\,\,\,- [ (\mathbf{D}(\mathcal{V}^{*(m)})\mathbf{n} \cdot \mathbf{n})\mathbf{n} - ((\mathbf{D}_m(v^{(m)})\mathbf{n}^{(m)} \cdot \mathbf{n}^{(m)})\mathbf{n}^{(m)}, \\
	\,\,\,\,\,\,\,\,\,\,\,\,\,\,\,\,\,\,\,\,\,\,\,\,\,\,\,\,\,\,\,\,\,\,\,\,\,\,\,\,\,\,- (\mathbf{D}_{m-1}(v^{(m-1)})\mathbf{n}^{(m-1)}\cdot \mathbf{n}^{(m-1)})\mathbf{n}^{(m-1)}) ]\} := d^{*(m)},\quad\text{on}\quad S_{F,T}, \\
	-\mathcal{Q}^{*(m+1)} + 2 \mathbf{D}(\mathcal{V}^{*(m+1)})\mathbf{n}\cdot \mathbf{n} = 2  [\mathbf{D}(\mathcal{V}^{*(m)})\mathbf{n}\cdot \mathbf{n} - ( \mathbf{D}_m(v^{(m)})\mathbf{n}^{(m)}\cdot \mathbf{n}^{(m)}, \\
	\,\,\,\,\,\,\,\,\,\,\,\,\,\,\,\,\,\,\,\,\,\,\,\,\,\,\,\,\,\,\,\,\,\,\,\,\,\,\,\,\,\,- \mathbf{D}_{m-1}(v^{*(m-1)})\mathbf{n}^{(m-1)}\cdot \mathbf{n}^{(m-1)})\mathbf{n}^{(m-1)}] := b^{*(m)} ,\quad\text{on}\quad S_{F,T}, \\
	\mathcal{B}_t^{*(m+1)} -  \triangle \mathcal{B}^{*(m+1)} =  (\triangle_m - \triangle_{m-1})B^{(m)} +  (\triangle_{m-1} - \triangle)\mathcal{B}^{*(m)} + \mathbf{B}\cdot(\nabla_m - \nabla_{m-1})\mathbf{v}, \\
	\,\,\,\,\,\,\,\,\,\,\,\,\,\,\,\,\,\,\,\,\,\,\,\,\,\,\,\,\,\,\,\,\,\,\,\,\,\,\,\,\,\,+ \mathbf{B}\cdot(\nabla_m - \nabla_{m-1})v^{*(m)} + \mathbf{B}\cdot\nabla_{m-1}\mathcal{V}^{*(m)} + B^{*(m)}\cdot(\nabla_m - \nabla_{m-1})\mathbf{v} + \mathcal{B}^{*(m)}\cdot\nabla_{m-1}\mathbf{v}, \\
	\,\,\,\,\,\,\,\,\,\,\,\,\,\,\,\,\,\,\,\,\,\,\,\,\,\,\,\,\,\,\,\,\,\,\,\,\,\,\,\,\,\,+ B^{*(m)}\cdot\nabla_m \mathcal{V}^{*(m)} + B^{*(m)}\cdot( \nabla_m - \nabla_{m-1} ) v^{*(m-1)} + \mathcal{B}^{*(m)}\cdot\nabla_{m-1} v^{*(m-1)} \\
	\,\,\,\,\,\,\,\,\,\,\,\,\,\,\,\,\,\,\,\,\,\,\,\,\,\,\,\,\,\,\,\,\,\,\,\,\,\,\,\,\,\, := g^{*(m)}  \quad \text{in}\quad Q_T, \\
	\mathcal{B}^{*(m+1)} = 0,\quad\text{on}\quad S_{F,T}, \\
	\mathcal{B}^{*(m+1)}(0) = 0,\quad\text{in}\quad\Omega\times\{t=0\}, \\
	\mathcal{V}^{*(m+1)} = \mathcal{B}^{*(m+1)} = 0,\quad\text{on}\quad S_{B}.
	\end{cases}
	\end{equation}
	
	Using Proposition~\ref{prop general stokes} and (\ref{Heat est}) with $v_0 = B_0 = 0$, we have
	\begin{equation} \label{Y ineq}
	\begin{split}
	Y_{m+1}^* &:= \|\mathcal{V}^{*(m+1)}\|_{H^{l+2, \frac{l}{2}+1}(Q_T)} + \|\mathcal{B}^{*(m+1)}\|_{H^{l+2, \frac{l}{2}+1}(Q_T)}  \\
	&\quad + \|\mathcal{Q}^{*(m+1)}\|^{(l)}_{Q_T} + \|\nabla \mathcal{Q}^{*(m+1)}\|^{(l)}_{Q_T} + \|\mathcal{Q}^{*(m+1)}\|_{W_2^{l+ \frac{1}{2}, \frac{l}{2}+ \frac{1}{4}}(S_{F,T})}  \\
	&\leq C_*(T) \{ \|f^{*(m)}\|^{(l)}_{Q_T} + \|g^{*(m)}\|^{(l)}_{Q_T} + \|\rho^{*(m)}\|_{W_2^{l+1, \frac{l+1}{2}}(S_{F,T})} + \|R^{*(m)}\|_{W_2^{0, \frac{l}{2}+1}(Q_T)} \\
	&\quad + T^{- \frac{l}{2}}\|R_t^{*(m)}\|_{L^2(Q_T)} + \|(b^{*(m)},d^{*(m)})\|_{W_2^{l+ \frac{1}{2}, \frac{l}{2}+ \frac{1}{4}}(S_{F,T})} + T^{- \frac{l}{2}}\|b^{*(m)}\|_{W_2^{ \frac{l}{2},0}(S_{F,T})} \},
	\end{split}
	\end{equation}
	where $C_{*}(T)$ is non-decreasingly time-dependent constant which means that it does not blow up as $T\rightarrow 0$. Now we estimate right hand side of (5.13). To estimate $\|f^{*(m)}\|^{(l)}_{Q_T}$, similar as (\ref{f_m}),
	\begin{equation} \label{f_m 1}
	\begin{split}
	\|   (\triangle_{m-1} - \triangle )\mathcal{V}^{*(m)} + (\triangle - \triangle_{m-1} )\mathcal{Q}^{*(m)} \|^{(l)}_{Q_T} \leq C(T,\delta_1,\delta_2,z^*)Y_m^*.
	\end{split}
	\end{equation}
	Using Lemma~\ref{larg lemma},
	\begin{equation} \label{f_m 2}
	\begin{split}
	\|   (\triangle_m - \triangle_{m-1})v^{(m)} + (\triangle_{m-1} - \triangle_m)q^{(m)} \|^{(l)}_{Q_T} \leq C(T,\delta_1,\delta_2,z^*)Y_m^*.
	\end{split}
	\end{equation}
	Similar as above two estimates,
	\begin{equation} \label{f_m 3}
	\begin{split}
	&\| \mathbf{B}\cdot(\nabla_m - \nabla_{m-1})\mathbf{B} + \mathbf{B}\cdot(\nabla_m - \nabla_{m-1})B^{*(m)} \\
	&+ B^{*(m)}\cdot(\nabla_m - \nabla_{m-1})\mathbf{B} + B^{*(m)}\cdot( \nabla_m - \nabla_{m-1} )B^{*(m-1)} \|_{Q_T}^{(l)} \\
	&\leq C(T,\delta_1,\delta_2,z^*)(Y_m^*)^2.
	\end{split}
	\end{equation}
	To control,
	\begin{equation} \label{four terms}
	\begin{split}
	\| \mathbf{B}\cdot\nabla_{m-1}\mathcal{B}^{*(m)} + {\mathcal{B}}^{*(m)}\cdot\nabla_{m-1}\mathbf{B} + B^{*(m)}\cdot\nabla_m \mathcal{B}^{*(m)} + \mathcal{B}^{*(m)}\cdot\nabla_{m-1} B^{*(m-1)} \|_{Q_T}^{(l)},
	\end{split}
	\end{equation}
	we estimate each four terms separately. We use Lemma~\ref{nonlinear est}, to get
	\begin{equation} \label{f_m 4}
	\begin{split}
	\| \mathbf{B}\cdot\nabla_{m-1}\mathcal{B}^{*(m)} \|_{Q_T}^{(l)} &\leq C ( T + T^{ \frac{1-l}{2}} + T^{ \frac{l}{2}} )\|\mathbf{B}\mathcal{G}^{(m)}\|_{H^{l+2, \frac{l}{2}+1}(Q_T)} \|\mathcal{B}^{*(m)}\|_{H^{l+2, \frac{l}{2}+1}(Q_T)} 
	\end{split}
	\end{equation}
	Meanwhile, on the RHS, we can use Lemma~\ref{fully nonlin est} to estimate,
	\begin{equation*}
	\begin{split}
	\|\mathbf{B}\mathcal{G}^{(m)}\|_{H^{l+2, \frac{l}{2}+1}(Q_T)} &\leq C \|\mathbf{B}\|_{H^{l+2, \frac{l}{2}+1}(Q_T)} \|\mathcal{G}^{(m)}\|_{H^{l+2, \frac{l}{2}+1}(Q_T)}  \\
	&\leq C \|\mathbf{B}\|_{H^{l+2, \frac{l}{2}+1}(Q_T)} \left( \|\mathcal{G}^{(m)} - \mathbf{I}\|_{H^{l+2, \frac{l}{2}+1}(Q_T)} + \|\mathbf{I}\|_{H^{l+2, \frac{l}{2}+1}(Q_T)} \right) .  \\
	\end{split}
	\end{equation*}
	Then, by Lemma~\ref{larg lemma}, we get estimate for the first terms in (\ref{four terms}).
	\begin{equation} \label{first term}
	\begin{split}
	\| \mathbf{B}\cdot\nabla_{m-1}\mathcal{B}^{*(m)} \|_{Q_T}^{(l)} \leq C(T,\delta_1,\delta_2,z^*) Y_m^* 
	\end{split}
	\end{equation}
	Other three terms can be controlled in similar ways as above and we get,
	\begin{equation} \label{other three}
	\begin{split}
	\| {\mathcal{B}}^{*(m)}\cdot\nabla_{m-1}\mathbf{B} \|_{Q_T}^{(l)} &\leq C(T,\delta_1,\delta_2,z^*) Y_m^* ,  \\
	\| B^{*(m)}\cdot\nabla_m \mathcal{B}^{*(m)} \|_{Q_T}^{(l)} &\leq C(T,\delta_1,\delta_2,z^*) Y_m^* ,  \\
	\| \mathcal{B}^{*(m)}\cdot\nabla_{m-1} B^{*(m-1)} \|_{Q_T}^{(l)} &\leq C(T,\delta_1,\delta_2,z^*) (Y_m^*)^2 .
	\end{split}
	\end{equation}
	Using (\ref{first term}), (\ref{other three}), (\ref{f_m 1}), (\ref{f_m 2}), and (\ref{f_m 3}), 
	\begin{equation} \label{f_m *}
	\|f^{*(m)}\|^{(l)}_{Q_T} \leq C(T,\delta_1,\delta_2,z^*) \left( Y_m^* + (Y_m^*)^2 \right) .
	\end{equation}
	Estimate of ${\|g^{*(m)}\|^{(l)}_{Q_T}}$ is similar to estimate of ${\|f^{*(m)}\|^{(l)}_{Q_T}}$ and we get,
	\begin{equation} \label{g_m *}
	\|g^{*(m)}\|^{(l)}_{Q_T} \leq C(T,\delta_1,\delta_2,z^*) \left( Y_m^* + (Y_m^*)^2 \right) .
	\end{equation}
	
	\noindent Lemma~\ref{larg lemma}, \ref{normal lemma}, and \ref{interpolation lemma} gives,
	\begin{equation} \label{rho_m *}
	\|\rho^{*(m)}\|_{W_2^{l+1, \frac{l+1}{2}}(S_{F,T})} \leq C(T,\delta_1,\delta_2,z^*)Y_m^*.
	\end{equation}
	
	\begin{equation} \label{R_m *}
	\|R^{*(m)}\|_{W_2^{0, \frac{l}{2}+1}(Q_T)} \leq C(T,\delta_1,\delta_2,z^*)Y_m^* +  C(T,\delta_1,\delta_2,z^*)(\varepsilon + C_{\varepsilon}T^{ \frac{1}{2}})Y_m^*.
	\end{equation}
	
	\begin{equation} \label{R_t_m *}
	T^{- \frac{l}{2}}\|R_t^{*(m)}\|_{L^2(Q_T)} \leq C(T,\delta_1,\delta_2,z^*)T^{ \frac{1-l}{2}}Y_m^*.
	\end{equation}
	
	\begin{equation} \label{bd_m *}
	\begin{split}
	&\|b^{*(m)}\|_{W_2^{l+ \frac{1}{2}, \frac{l}{2}+ \frac{1}{4}}(S_{F,T})} + \|d^{*(m)}\|_{W_2^{l+ \frac{1}{2}, \frac{l}{2}+ \frac{1}{4}}(S_{F,T})} + T^{- \frac{l}{2}}\|b^{*(m)}\|_{W_2^{ \frac{l}{2},0}(S_{F,T})} \\
	&\leq C(T,\delta_1,\delta_2,z^*)T^{ \frac{1-l}{2}}Y_m^*,
	\end{split}
	\end{equation}
	Now we put (\ref{f_m *}) -- (\ref{bd_m *}) together to derive, 
	\begin{equation}
	Y_{m+1}^* \leq \chi Y_m^*,
	\end{equation}
	where $\chi < 1$ if we pick a $\varepsilon$ and sufficiently small $T$ which is smaller than $T_1$ of (\ref{delta 12}). Hence, by contraction mapping principle, we solve (\ref{fully nonlin system}). So far, we proved Theorem \ref{main theorem}, except $\nabla\cdot H = 0$. 
	
	\subsection{Divergence free of $H$}
	System (\ref{fully nonlin system}) (which was solved in section 5) does not necessarily satisfy divergence-free condition of magnetic field $H$. We recover our system in Eulerian coordinate and claim that divergence-free property of $H(t)$ propagates from initial condition $\nabla\cdot H_0 = 0$. We appeal to maximum principle of convection-diffusion equation.
	\begin{equation*}
	\begin{cases}
	H_t + (u\cdot\nabla)H - (H\cdot\nabla)u =  \triangle H, \\
	\nabla\cdot u = 0 ,  \\
	\nabla\cdot H_0 = 0.
	\end{cases}
	\end{equation*}
	Taking divergence to above equation and using notation  $\mathcal{H} := \nabla\cdot H $,
	$$
	\mathcal{H}_t + (u\cdot\nabla)\mathcal{H} + (\nabla u):(\nabla H)^t - (\nabla H):(\nabla u)^t - (H\cdot\nabla)(\nabla\cdot u) =  \triangle \mathcal{H},
	$$
	where $A:B := \sum_{i,j} A_{ij} B_{ij}$. Hence,
	$$
	\mathcal{H}_t + (u\cdot\nabla)\mathcal{H} -  \triangle \mathcal{H} = 0.
	$$
	Then by maximum principle of convection-diffusion equation, 
	$$
	\|(\nabla\cdot H) (t)\|_{L^\infty} = \|\mathcal{H}(t)\|_{L^\infty} \leq \|\mathcal{H}(0)\|_{L^\infty} = \|\nabla\cdot H_0\|_{L^\infty} = 0,
	$$
	during the time interval for the solution $H$.

	\section*{Acknowledgments}
	The author was partially supported by NSF grant DMS-1211806.

	\vspace{1cm}
	\indent\indent \author{$\textsc{Donghyun\,\,\,Lee}$}\\
	\indent \address{$\textsc{University\,\,\,of\,\,\,Wisconsin-Madison,\,\,\,Van Vleck Hall, 480 Lincoln Drive, Madison, Wi  53706}$} \\
	\indent \email{$\mathtt{dlee374@wisc.edu}$}

\end{document}